\documentclass[12pt]{article}
\usepackage{epsfig}
 \usepackage{graphics}
\usepackage{amsmath,amssymb}
\usepackage{subfigure}
\usepackage{graphicx}

\newcommand{\thefont}[2]{\fontsize{#1}{#2}\fontshape{n}\selectfont}
\newcommand{\1}{\rlap{\thefont{10pt}{12pt}1}\kern.16em\rlap{\thefont{11pt}{13.2pt}1}\kern.4em}
\newcommand{\RR}{\ensuremath{\mathbb{R}} }

\newcommand{\NN}{\ensuremath{\mathbb{N}} }

\newcommand{\bbeta}{\ensuremath{\boldsymbol{\beta}} }

\newcommand{\balpha}{\ensuremath{\boldsymbol{\alpha}} }

\newcommand{\bx}{\ensuremath{\mathbf{x}} }

\newcommand{\bF}{\ensuremath{\mathbf{f}} }
\newcommand{\bY}{\ensuremath{\mathbf{Y}} }

\newcommand{\PR}{\ensuremath{{\mathbb P}}}
\newcommand{\HH}{\ensuremath{{\mathcal H}}}
\newcommand{\FF}{\ensuremath{{\mathcal F}}}
\newcommand{\EE}{\ensuremath{{\mathbb E}}}
\newcommand{\opO}{\ensuremath{\mathcal O}}

\newtheorem{theo}{Theorem}
\newtheorem{prop}{Proposition}
\newtheorem{lemma}{Lemma}
\newtheorem{definition}{Definition}

\newtheorem{coro}{Corollary}

\numberwithin{equation}{section}
\numberwithin{ass}{section}
\numberwithin{theo}{section}
\numberwithin{prop}{section}
\numberwithin{lemma}{section}
\numberwithin{definition}{section}
\numberwithin{rmq}{section}

\title{Smoothing under diffeomorphic constraints with homeomorphic splines}

\author{{\em J\'er\'emie Bigot \& S\' ebastien Gadat} \\
Institut de Math\'ematiques de Toulouse \\
 Universit\'e Paul Sabatier \\
 F-31062 Toulouse Cedex 9, France \\
{\tt email: \{Jeremie.Bigot,Sebastien.Gadat\}@math.univ-toulouse.fr}
}

\begin{document}

\thispagestyle{empty}

\maketitle

\begin{abstract}
In this paper we introduce a new class of diffeomorphic smoothers based on general spline smoothing techniques and on the use of some tools that have been recently developed in the context of image warping to  compute smooth diffeomorphisms. This diffeomorphic spline  is defined as the solution of an ordinary differential equation governed by an appropriate time-dependent vector field. This solution has a closed form expression which can be computed using classical unconstrained spline smoothing techniques. This method does not require the use of quadratic or linear programming under inequality constraints and has therefore a low computational cost. In a one dimensional setting incorporating diffeomorphic constraints is equivalent to impose monotonicity. Thus, as an illustration, it is shown that such a monotone spline can be used to monotonize any unconstrained estimator of a regression function, and that this monotone smoother inherits the convergence properties of the unconstrained estimator. Some numerical experiments are proposed to illustrate its  finite sample performances, and to compare them with another monotone estimator. We also provide a two-dimensional application on the computation of diffeomorphisms for landmark and image matching.
\medskip


\end{abstract}

\section{Introduction}

Spline smoothing is widely used in many different areas to study the relationship
between a response variable $Y$ and an independent variable $X$, see e.g. Wahba \cite{wah} for a detailed presentation, and have many applications in approximation problems, see e.g. Duchon \cite{duchon}, Quak \& Schumaker \cite{quak}, de Boor \& Schumaker \cite{SchumdeBoor} or Lopez de Silanez \& Apprato \cite{MCL} . In many fields of interest, including physical and medical sciences, one is often interested in imposing a monotonic relationship between two such variables. Typical examples include the analysis of dose-response curves in pharmakinetics (Kelly \& Rice~\cite{kel}), growth curves in biology and many specific practical problems discussed in the literature cited below.  Note that without loss of generality, monotone smoothing is  considered in this paper  as the problem of computing an increasing function. For calculating a decreasing smoother one can simply reverse the ``$X$ axis'' and then apply the same methodology.

Unconstrained spline smoothing consists in minimizing over an appropriate functional space $\mathcal{F}$ a criterion that represents a balance between two terms : fidelity to the data and smoothness of the fitting spline. If $\mathcal{F}$ is a Reproducing Kernel Hilbert Spaces (RKHS), then it is well known (see e.g.  Wahba \cite{wah}) that a closed form solution can be computed by solving a simple linear system of equations. The simplest idea that comes to mind to incorporate monotonicity constraints is to restrict the search space $\mathcal{F}$ to a subset of monotone functions, and then to take as a monotone smoother the function which minimizes the same criterion over this restricted space. Existence, approximation properties and the numerical computation of such monotone smoothers have been widely studied, see e.g. Utreras \cite{utre}, Andersson \& Elfving  \cite{and}, Elfving \& Andersson \cite{elf}, Irvine, Marin \& Smith \cite{irv}, Beatson \& Ziegler \cite{beat}.

However, finding the exact solution of a smoothing problem in a constrained space is generally a difficult task and most existing algorithms only produce approximate solutions. As a closed form solution for such monotone smoothers does not exist in general, their numerical computation is generally done by determining the fitted values of the  smoothing spline on a finite set of points (usually the observed covariates) and uses a set of inequality constraints to impose restrictions on the value of the fitted function at these points. However, the algorithms used to compute these estimators can be computationally intensive since they involve a large set of inequality  constraints (see e.g.  Schwetlick \& Kunert \cite{schwe}, He \& Shi \cite{he}, Turlach  \cite{turl} and the discussion therein).

In this paper, we introduce a new class of monotone smoothers that have a closed form expression which depends on the underlying RKHS, and which can be computed using classical unconstrained spline smoothing techniques. Thus, unlike some monotone smoothers, this method does not require the use of quadratic programming under linear constraints and has therefore a low computational cost. Our approach is based on tools that have been recently developed in the context of image warping for the construction of diffeomorphisms in two or three dimensions  (see e.g. Trouv\'e \cite{trou},  Miller, Trouv\'e  \& Younes \cite{miltrou}, Glaun\`es \cite{glau},  Younes \cite{youbook}, \cite{you}, see also Apprato \& Gout \cite{app} for the use of diffeomorphism in spline approximation).  Trouv\'e, Younes and their collaborators have proposed to compute smooth diffeomorphisms as the  solutions of a new class of ordinary differential equation governed by a time-dependent velocity field. In a one-dimensional (1D) setting, it is easy to see that a diffeomorphism is a smooth and monotone function, and thus the main idea of this paper is to adapt such tools for the construction of 1D monotone smoothers. Our approach also yields a new method to compute smooth diffeomorphisms for the alignment of landmarks in a 2D or 3D setting. Some examples in the 2D case for image warping (see e.g. Bigot {\it et al.} \cite{bigot}) are given in the section on numerical experiments, but for simplicity our theoretical results are presented in a 1D setting.

Our main contributions are the following: first we show how one can generate a strictly monotone function as the solution of an ordinary differential equation (ODE). We also prove that for some functional classes, any monotone function can be represented as the solution of such an ODE. Secondly, a new criterion to fit a monotone spline is proposed and we show that the minimizer of this criterion has a simple closed form expression. As an illustrative example, we explain how the overall methodology can be applied to the problem of monotonizing any estimator of a regression function. Indeed, in statistics, a possible smoothing method under shape constraints consists in first using an unconstrained estimator (such as a spline, wavelets or kernel smoother) and then projecting the resulting curve estimate onto a constrained subspace of regression functions  which is usually a convex set (see e.g. Mammen, Marron, Turlach \& Wand~\cite{mam1} and Mammen \& Thomas-Agnan~\cite{mam2}). For the problem of monotone regression, this approach is generally referred to  as  smooth and then monotonize. However, as pointed out by Gijbels~\cite{gib} many of these monotone estimates appear less smooth than the unconstrained estimates due to the projection step. Moreover, it is not clear how one can compute numerically the projection of a curve estimate onto a constrained subspace for any unconstrained estimator. Our monotone estimator does not suffer from these two drawbacks since it can be easily computed, and it yields surprisingly very smooth estimates.

The remainder of the paper is structured as follows:  section \ref{secRKHS} gives a brief overview of RKHS  and the general spline smoothing problem. In section \ref{secode}, we show how one can generate a strictly monotone function as the solution of an ODE. In section \ref{sechomeo}, we propose a new class of monotone smoothers that we call homeomorphic smoothing splines.  In Section \ref{sec:stat}, we apply this methodology to nonparametric regression under monotonicity constraints.  Section  \ref{secsimu} presents  a short Monte Carlo study on the efficiency of this approach and a comparison with another constrained estimate. Finally, another application of homeomorphic splines is presented in a 2D setting for matching problems involving the alignment of landmarks. The Appendix provides the proofs of the main results.

\section{The general spline smoothing problem in  the 1D case} \label{secRKHS}

Let $\HH_{K}$ a be RKHS of functions in
$L^{2}(\RR)$ with positive definite kernel $K$, meaning that for all $x \in \RR$ there exists an element $K_{x}(.)$ such that
for all $g$ in $\HH_{K}$, we have $g(x) = \langle K_{x}(.),g\rangle_{K} $ where
$\langle .,.\rangle_{K}$ denotes the scalar product on $\HH_{K}$ whose derived
norm is $\|. \|_{K}$ (for more details on RKHS we refer to Atteia \cite{Atteia1}, Atteia \& Gaches \cite{Atteia2}, Aronszajn \cite{Aron}, Duchon \cite{duchon}, Wahba \cite{wah},  Berlinet \& Thomas-Agnan \cite{bertho}).  Let $\psi_{1},\ldots,\psi_{M}$ be functions (not necessarily in $L^{2}(\RR)$)
such that for any set of distinct points $x_{1},\ldots,x_{n}$ in $\RR$, the
matrix $T$ with elements $T_{i,j} = \psi_{j}(x_{j})$ has full rank $M < n$. Let
$\tilde{\HH} = Span \{ \psi_{j} \}_{j=1,\ldots,M} + \HH_{K}$ and assume that we have $n$ distinct pairs of points $(x_{i},y_{i}) \in \RR^{2}$.  Then the general spline smoothing problem is to find the minimizer
\begin{equation}
\tilde{h}_{n,\lambda} = \arg \min_{\tilde{h} \in \tilde{\HH}} \; \frac{1}{n} \sum_{i=1}^{n} (\tilde{h}(x_{i})-y_{i})^{2} + \lambda \|h \|^{2}_{K}, \label{eqsmoothsplinegen}
\end{equation}
for any $\tilde{h} \in \tilde{\HH}$ of the form $\tilde{h}(x) = \sum_{j=1}^{M} a_{j}  \psi_{j}(x) +h(x)$,
where $h \in \HH_{K}$.
It is well known that the solution of this smoothing problem is unique and of the form: $\forall x \in \RR, \;  \tilde{h}_{n,\lambda}(x) = \sum_{j=1}^{M} \alpha_{j}
\psi_{j}(x) + \sum_{i=1}^{n} \beta_{i} K(x,x_{i}),
$ where $\alpha,\beta$ are solutions of a simple linear system of
equations (see e.g. Wahba \cite{wah}).

Throughout this paper, we will assume $M = 2$ with $\psi_{1}(x) = 1$ and
$\psi_{2}(x) = x$, these conditions are convenient for satisfying the uniform Lipschitz condition stated in Lemma \ref{lemLipBound} (see the Appendix).

As an example of a RKHS, we will often use the Sobolev space of order $m \in \NN^{\ast}$ (see  Berlinet and Thomas-Agnan \cite{bertho})
$\HH_{K} = \HH^{m}(\RR)$ endowed with the norm
$ \|h\|^2_{\HH^{m}} = \int_{-\infty}^{+\infty} |h(x)|^{2} dx +   \int_{-\infty}^{+\infty} |h^{(m)}(x)|^{2} dx.$ With this choice for the $\psi_{j}$'s,  any function $\tilde{h} \in
\tilde{\HH}$ is of the form $\tilde{h}(x) = a_{1} + a_{2}x + h(x)$ where $a_{1},a_{2}
\in \RR$ and $h \in \HH_{K}$.  Hence, one can define a norm in $\tilde{\HH}$ by
setting $$ \| \tilde{h}\|_{\tilde{\HH}} = \max(|a_{1}|,|a_{2}|) +    \|h\|_{K}.$$ Note that $\tilde{\HH}$ is a Banach space for this norm. From now, to simplify
the notations we will omit the superscript $\tilde{\HH}$ and write $\|
\tilde{h}\|_{\tilde{\HH}} = \| \tilde{h}\|$.


\section{Differential Equation to generate monotone functions} \label{secode}

The smoothing spline $\tilde{h}_{n,\lambda}$ defined previously is
not necessarily a monotone function. We thus propose to use a connection between
monotone functions and time-dependent vector fields to incorporate monotonicity
constraints into the computation of $\tilde{h}_{n,\lambda}$.

\subsection{Generating monotone functions}

Let us explain the basic ideas (as described e.g. in Younes \cite{youbook})
to generate smooth diffeomorphisms. Take any $v \in \mathcal{C}^1(\RR,\RR)$ with $\|v'\|_{\infty} < +\infty$. Then if $\epsilon>0$ is chosen sufficiently
small the perturbation $\phi = Id+\epsilon v$ of the identity function, is a
strictly increasing \textit{small} diffeomorphism. Now, if $v_{t_1},\ldots,v_{t_p}$
are continuously differentiable functions on $\RR$ and  if $\epsilon > 0$ is
sufficiently small such that $Id+ \epsilon v_{t_k}$ are small diffeomorphisms on $\RR$, then
we can construct the following sequence of diffeomorphisms $\phi_{t_p} =
(Id+\epsilon v_{t_{p-1}}) \circ \ldots \circ  (Id +\epsilon v_{t_1})$. Then, note
that $\phi_{t_{p+1}} =  (Id +\epsilon v_{t_p}) \circ \phi_{t_p} =  \phi_{t_p} + \epsilon v_{t_p} \circ \phi_{t_p}$ which can also be written as
\begin{equation}
\forall x \in \RR \qquad \frac{\phi_{t_{p+1}}(x) -  \phi_{t_p}(x)}{\epsilon} = v_{t_p}
\left[\phi_{t_p}(x)\right]. \label{eqodedis}
\end{equation}
As $\epsilon \to 0$ and $\lim_{p \to + \infty} t_{p+1}-t_{p}=0$ , (\ref{eqodedis}) looks like a discretized version of an
inhomogeneous differential equation of the following form (by introducing  a
continuous time variable $t$):
\begin{equation}
\label{edo} \frac{\text{d} \phi_{t}}{\text{d} t} = v_{t}( \phi_{t}).
\end{equation}

In the sequel, $v$ will be a function of two variables $(t,x)$, and for a fixed $t$ we will use the notation  $x \mapsto v_t(x)$ to refer to the application $x \mapsto v(t,x)$. The variable $t$ varies in the finite time interval $[0;1]$ while $x$ belongs to $\mathbb{R}$. Similarly $\phi$ depends both on the time $t$ and the variable $x$, and $\phi_t(x)$ will refer to $\phi(t,x)$. Thus, equation (\ref{edo}) is equivalent to:
$
\forall x \in \mathbb{R},  \frac{\text{d} \phi(t,x)}{\text{d} t} = v(t, \phi(t,x)).
$

As we will see, under mild conditions on the time-dependent vector field $(v_{t})_{t \in [0,1]}$, the solution of the above ODE is a diffeomorphism at all time $t$ and thus a monotone function. The main idea of this paper is thus the following: we transfer the problem of computing a monotone spline from a set of $n$ data points $(x_i,y_i)_{i = 1 \dots n}$ to the problem of computing an appropriate vector field $(v_t^{n})_{t \in [0,1]}$ which depends on these data points. A monotone smoother $f_{n}$ is then defined as the solution at time $t=1$ of the ODE  (\ref{edo}) governed by the vector field  $(v_t^{n})_{t \in [0,1]}$, i.e $f_n(x) = \phi_1(x)$ with $\phi_0(x) = x$. The main advantage of this approach is that the computation of $(v_t^{n})_{t \in [0,1]}$ will be obtained from an unconstrained smoothing problem, and therefore the calculation of $f_{n}$ only requires to run an ODE without imposing any specific constraints. This yields a fitting function $f_{n}$ which is guaranteed to be monotone.

\subsection{Vector fields and ODE}

Following the notations in Younes \cite{you}, let us state several definitions.\\

\begin{definition}[$\mathcal{X}^{1},\mathcal{X}^{2}$ and $\mathcal{X}$]
$\mathcal{X}^{1}$  is the space of time-dependent vector fields $(v_{t} \in \tilde{\HH} ,t \in
[0,1])$ such that $ \| v \|_{\mathcal{X}^{1}} =_{def} \int_{0}^{1} \|v_{t}\|dt < + \infty$.  $\mathcal{X}^{2}$ is the space of time-dependent vector fields $(v_{t} \in \tilde{\HH} ,t \in
[0,1])$ such that $\| v \|_{\mathcal{X}^{2}} =_{def} \int_{0}^{1} \|v_{t}\|^{2}dt < + \infty.$ Finally, $\mathcal{X}$ is the set of all time-dependent vector field $(v_{t} \in
\tilde{\HH} ,t \in [0,1])$.\\
\end{definition}

The definitions of $\mathcal{X},\mathcal{X}^{1},\mathcal{X}^{2}$ are classical in the context of time-dependent PDE which are formulated as Banach space-valued functions, see e.g. Renardy \& Rogers \cite{rr}. Note that by the Cauchy-Schwarz inequality, $\| v \|_{\mathcal{X}^{1}}^{2} \leq \|
v \|_{\mathcal{X}^{2}}$, and thus $ \mathcal{X}^{2} \subset \mathcal{X}^{1} \subset
\mathcal{X}$.  For  $v \in \mathcal{X}^{1}$, we formally define an ODE governed
by the time-dependent vector field $(v_{t},t \in [0,1])$ as
\begin{equation}
\frac{\text{d} \phi_{t}}{\text{d} t} = v_{t}( \phi_{t}). \label{eqodecont}
\end{equation}

\begin{definition}
Let $\Omega = [0,1]$. A function $t \mapsto \phi_{t}$ is called a solution of
the equation  (\ref{eqodecont}) with initial condition the identity  if for all $x \in \Omega$, $t \mapsto\phi_{t}(x)$ is a continuous function from $[0,1]$ to $\RR$, $\phi_{0}(x) = x$ for all $x \in \Omega$,
and  for all $t \in [0,1]$ and all $x \in \Omega$, $\phi_{t}(x) = x + \int_{0}^{t} v_{s}(\phi_{s}(x)) ds.$\\
\end{definition}

The following theorem, whose proof is deferred to the appendix, shows that the solution of the equation  (\ref{eqodecont})
is unique and is a homeomorphism for all time $t \in [0,1]$.\\

\begin{theo} \label{thode}
Assume that the kernel $K$ is bounded on $\RR^{2}$, and that there exists a constant $C_{1}
$ such that for any $h \in \HH_{K}$
\begin{equation}
|h(x) - h(y)| \leq C_{1} \|h\|_{K}  |x-y|. \label{eqassLip}
\end{equation}
Let $v \in \mathcal{X}^{1}$. Then, for all $x \in \Omega$ and $t \in [0,1]$, there exists a unique
solution of (\ref{eqodecont}) with initial
 condition  the identity. Besides, for all $t \in [0,1]$, $\phi_{t}$ is a
  homeomorphism from $\Omega$ to $\phi_{t}(\Omega)$.\\
\end{theo}

The above uniformly Lipschitz assumption (\ref{eqassLip}) for the kernel $K$ is not restrictive as it is satisfied in many cases of interest. Indeed, observe that for for any $h \in \HH_{K} $:
$$
|h(x)-h(y)| = | \langle K(.,x)-K(.,y),h \rangle_{\HH_{K}}| \leq \|h\|_{K} \|K(.,x)-K(.,y) \|_{K}.
$$
If $K$ is a radial kernel of the form $K(x,y)=k(|x-y|)$ for some function $k : \RR \to \RR$, then
the above equation implies that $|h(x)-h(y)| \leq 2 \|h\|_{K}  |k(0) - k(|x-y|)|$. Hence, $h$ satisfies equation (\ref{eqassLip}) provided $k$ is uniformly Lipschitz on $\RR$. This is the case for a Gaussian kernel: $k(|x-y|) = e^{-|x-y|^2/2\sigma^2}$, and also in the Sobolev case  where $K$ is given by (see e.g. \cite{bertho}):
$
K(x,y) = k_m(|x-y|) := \sum_{k=0}^{m-1} \frac{\exp \left[-|x-y| \exp(i\frac{\pi}{2m}+\frac{k \pi }{m} - \frac{\pi}{2})\right]}{2 m \exp((2m-1) (i \frac{\pi}{2m}+\frac{k \pi }{m})}.
$ \\

\noindent {\bf Remark:} this framework can be extended to a 2D setting for generating diffeomorphism of $\RR^{2}$. For this, let $ \tilde{\HH}_{2}$ denote a set of smooth functions from $\RR^{2}$ to $\RR^{2}$ (see Section \ref{sec:H2} for an example) and  define $ \tilde{\mathcal{X}}= \{ (v_{t},t \in [0,1]) \mbox{ with }  v_t \in \tilde{\HH}_{2} \mbox{ for all } t \in [0,1]  \}$. Let Let $\Omega$ be an open subset of $\RR^{2}$ and define an ODE governed by the time-dependent vector field $(v_{t},t \in [0,1]) \in  \tilde{\mathcal{X}}$ as:
$
\frac{\text{d} \phi^{v}_{t}}{\text{d} t} = v_{t}( \phi^{v}_{t}),
$
with $v_{0}(x) = x$. Using arguments in the proof of Theorem \ref{thode} and assuming that the functions $h$ in $\tilde{\HH}_{2}$ are sufficiently smooth and satisfy a uniform Lipschitz condition of the type (\ref{eqassLip}), then one can easily show that the solution of such an ODE is unique and is a diffeomorphism from $\Omega$ to $\phi_{t}(\Omega) \subset \RR^{2}$ for all time $t \in [0,1]$.

\section{Homeomorphic smoothing splines} \label{sechomeo}


\subsection{A connection between monotone functions and time-dependent vector fields}

A natural question is to ask if any monotone function can be written as the solution of an ODE governed by  a time-dependent vector field. First, consider the case where $\HH_{K} = \HH^{m}(\RR)$ and $f$ belongs to the Sobolev space
$$
H^{m}([0,1]) = \{f : [0,1] \to \RR, \text{$f^{(m-1)}$ is absolutely continuous with} \int_{0}^{1} |f^{(m)}(x)|^{2}dx < + \infty \}.
$$
Then, if $f$ is monotone, one of our main results is the following theorem which states that $f$ can be represented as the solution at time $t=1$ of an ODE:\\

\begin{theo} \label{theovf}
Assume that $\tilde{\HH} = Span\{1,x\} + \HH^{m}(\RR)$. Let $m \geq 2$ and $f  \in   H^{m}([0,1])$ be such that $f'(x) > 0$ for all $x \in [0,1]$ and define $\phi_{t}(x) = tf(x) + (1-t)x$, for all $t \in [0,1]$. Then, there exits a time-dependent vector field $(v_t^{f})_{t \in [0,1]}$ depending on $f$, such that  $v^{f}_{t} \in \HH^{m}(\RR)$ for all $t \in [0,1]$ and which satisfies
$\phi_1 = \phi_0 + \int_{0}^{1} v^{f}_{t}(\phi_{t}) dt$, and thus
$
f(x) = \phi_1(x) =  x + \int_{0}^{1} v^{f}_{t}(\phi_{t}(x)) dt.
$
Moreover for all $t \in [0,1]$ one has that
\begin{equation}
v^{f}_{t}( \phi_{t}(x) ) = v^{f}_{t}(   tf(x) + (1-t)x ) = f(x) - x = \frac{\text{d} \phi_{t}}{\text{d} t}(x) \mbox{ for all } x \in [0,1]. \label{eqvf1}
\end{equation}
For all $t \in [0,1]$, the function $v^{f}_{t}$ can be chosen as the unique element of minimal norm in $H^{m}(\RR)$ which satisfies equation (\ref{eqvf1}).\\
\end{theo}

To the best of our knowledge, this representation of a monotone function by such an ODE has not been used before.  The formulation (\ref{eqvf1}) suggests the following trick to compute a monotone smoother from a set of $n$ data points $(x_{i},y_{i}) \in [0,1] \times \RR $: if one considers $y_{i}$ as an approximation of $f(x_{i})$ for some function $f$, then to obtain a good approximation of $v^{f}$, one can  use the  $y_{i}$ 's to compute a vector field   $v^{n}$ that satisfies roughly the interpolating conditions (\ref{eqvf1}) at the design points. More precisely, at any time $t$, the vector field $v^{n}_{t}$ is obtained by smoothing the ``data''  $(t y_{i} + (1-t)x_{i},y_{i} - x_{i}),
\; i=1,\ldots,n$. Finally, to compute a monotone smoother $f_{n}$ we just have to run the ODE (\ref{eqodecont}) with the vector field $v^{n}$.

When $\HH_{K} \neq \HH^{m}(\RR)$, it is not clear if one can obtain a general correspondence between  monotone functions and their representation via a vector field $v \in \mathcal{X}^{2}$. However, we believe that the proof of Theorem \ref{theovf} could be adapted to other RKHS.

\subsection{A new monotone smoothing spline}

Let $(x_{i},y_{i}), i =1,\ldots,n$ be a set of data points with $x_{i}  \in [0,1]$ and $y_{i} \in  \RR$. A new smoothing spline problem under monotonicity constraints can be formulated in the following way: for a time-dependent vector field $v \in \mathcal{X}$, define the ``energy''
\begin{equation}
E_{\lambda}(v) =  \int_{0}^{1} \frac{1}{n} \sum_{i=1}^{n}\left(y_{i}-x_{i} -v_{t}(ty_{i}+(1-t)x_{i}) \right)^{2}dt + \lambda \int_{0}^{1} \| h_{t} \|^{2}_{K} dt,\label{eqhomeospline1}
\end{equation}
where $v_{t}(x) = a_{1}^{t} + a_{2}^{t}x + h_{t}(x)$, and $\lambda > 0$ is a regularization parameter. Then, take $ v^{n,\lambda} = \arg \min_{v \in \mathcal{X}} E_{\lambda}(v),$ and a monotone smoother $f_{n,\lambda}^{c}$  is obtained  by taking $f_{n,\lambda}^{c}(x) = \phi^{v^{n,\lambda}}_{1}(x) = x + \int_{0}^{1} v^{n,\lambda}_{t}(\phi^{v^{n,\lambda}}_{t}(x)) dt.$ The following proposition gives sufficient conditions for the existence of $v^{n,\lambda} $:\\

\begin{prop}  \label{propexist}
Assume that the conditions of Theorem \ref{thode} are satisfied. Assume that $n >  2$ and that the kernel $K: \RR^{2} \to \RR$ is continuous. Suppose that the $x_{i}$'s and the $y_{i}$'s are such that the $n$ ``design points'' $ty_{i}+(1-t)x_{i}$ are distinct in $\RR$ for any $t \in [0,1]$. Then, the optimization problem (\ref{eqhomeospline1}) has a unique solution $v^{n,\lambda} \in \mathcal{X}$ such that at each time $t \in [0,1]$, $v_{t}^{n,\lambda}$ is the solution of the following standard unconstrained smoothing problem: find $v_{t} \in   \tilde{\HH}$ which minimizes
\begin{equation}
E^{t}_{\lambda}(v_{t}) = \frac{1}{n}\sum_{i=1}^{n} \left(y_{i}-x_{i} -v_{t}(ty_{i}+(1-t)x_{i}) \right)^{2} + \lambda  \|h_{t} \|^{2}_{K}, \label{eqener0}
\end{equation}
where $v_{t}(x) = a_{1}^{t} + a_{2}^{t}x + h_{t}(x)$. Moreover $v^{n,\lambda} \in \mathcal{X}^{2}$, and $f_{n,\lambda}^{c}$ is a monotone function on $[0,1]$.\\
\end{prop}

Let us remark that if one defines $\tilde{\HH}^{c}$ as the subspace of functions  $f \in \tilde{\HH}$ such that $f$ is a strictly monotone function on $[0,1]$, then a monotone smoother $\tilde{f}_{n,\lambda}^{c}$ can de defined by mininizing the classical spline smoothing criterion over the restricted space $\tilde{\HH}^{c}$ i.e.
\begin{equation}
\tilde{f}_{n,\lambda}^{c} = \arg \min_{\tilde{h} \in \tilde{\HH}^{c}} \; \frac{1}{n} \sum_{i=1}^{n} (\tilde{h}(x_{i})-y_{i})^{2} + \lambda \|h \|^{2}_{K}. \label{eqsmoothsplineconst}
\end{equation}
General smoothing splines problems under under shape constraints such as monotonicity have been studied in detail in Utreras \cite{utre}.  Theorems proving the existence, uniqueness and general results concerning the characterization of $\tilde{f}_{n,\lambda}^{c}$ are given in  Utreras \cite{utre}, together with a study of the convergence rate of $\tilde{f}_{n,\lambda}^{c}$  in a nonparametric regression setting. Hence, it would be interesting to study the relationship that may exist between the estimators $\tilde{f}_{n,\lambda}^{c}$ and $f_{n,\lambda}^{c}$. However, we believe that this problem is not an easy task which is beyond the scope of this paper.


\subsection{Computational aspects and the choice of $\lambda_{n}$} \label{seccomp}
\paragraph{Numerical computation}
The optimization problem (\ref{eqener0})  amounts to solve, at each time $t \in [0,1]$, a simple finite-dimensional least-square problem which yields a very simple algorithm to compute a smooth increasing function: choose a discretization $t_{k} = \frac{k}{T},k=0,\ldots,T-1$ of the time-interval $[0,1]$ (we took $T=30$) and set $\phi^{ t_{0}}_{n,\lambda}(x) = x$ for $x \in [0,1]$. Then
repeat for $k = 0,\ldots,T-1$: find the solution $v_{t_{k}}^{n,\lambda}$ of the unconstrained smoothing problem (\ref{eqener0}) for each $t = t_{k}$, and then compute $\phi^{ t_{k+1}}_{n,\lambda}(x) = \phi^{ t_{k}}_{n,\lambda}(x) +\frac{1}{T} v_{t_{k}}^{n,\lambda}\left(\phi^{t_{k}}_{n,\lambda}(x)\right)$. The proposed numerical scheme is based on
\begin{equation}
\label{eq:euler}\phi^{ t_{k+1}}_{n,\lambda} = (Id+  \frac{1}{T} v_{t_{k}}^{n,\lambda}) (\phi^{ t_{k}}_{n,\lambda})
\end{equation}
which replaces the theoretical relation $\phi^{ t_{k+1}}_{n,\lambda} = \phi^{ t_{k}}_{n,\lambda}  + \int_{t_k}^{t_{k+1}} v_{u}^{n,\lambda} (\phi^{u}_{n,\lambda}) du$. Remark that equation (\ref{eq:euler}) shows that if  $ \frac{\|v^{n,\lambda}_t\|}{T} < 1$ for all $t$, then $\phi_{n,\lambda}^{t_{k+1}}$ remains monotone provided $\phi_{n,\lambda}^{t_k}$ is monotone. This condition is not really restrictive since we have shown that $(v_{t}^{n,\lambda})_{t \in [0,1]}$ is in $\mathcal{X}^{2}$ and that $t \to v_{t}^{n,\lambda}$ is a continuous map on $[0;1]$. Thus, our estimator based on the Euler scheme (\ref{eq:euler}) remains monotone if $T$ is chosen sufficiently large, namely greater than $\sup_{t \in [0,1] } \|v^{n,\lambda}_t\|$.

Another important question is the error made using the Euler discretization scheme instead of the correct ODE This point is left open since it is far beyond the scope of this paper but the use of the Gronwall Lemma should enable to derive upper bound between the theoretical $\phi^1_{n,\lambda}$ and the approximated one derived from (\ref{eq:euler}).
\paragraph{Choice of the regularization parameter}
A fundamental issue is the choice of the regularization parameter $\lambda$.  In our simulations, we have obtained good results via an empirical choice of $\lambda$ inspired by the generalized cross-validation (GCV) criterion of Craven \& Wahba \cite{cra}, see also Girard \cite{girard} for fast cross-validation methods. For $i = 1,\ldots,n$ and $t \in [0,1]$, define $\hat{X}_{i}^{t} = t y_{i}+(1-t)x_{i}$ and $\hat{Y}_{i} =y_{i}-x_{i} $. Then, note that at each time $t \in [0,1]$ the smoothing spline $v_{t}^{n,\lambda}$ evaluated at the ``design points'' $\hat{X}_{i}^{t},\ldots,\hat{X}_{n}^{t}$ is a linear function of the observations $\hat{Y}_{1},\ldots,\hat{Y}_{n}$, i.e. there exists a matrix $A_{\lambda,t}$ such that $\mathbf{v}_{t}^{n,\lambda} = \left(v_{t}^{n,\lambda}(\hat{X}_{1}^{t}),\ldots,v_{t}^{n,\lambda}(\hat{X}_{n}^{t})\right)' = A_{\lambda,t} \bY$, with $\bY = (\hat{Y}_{1},\ldots,\hat{Y}_{n})'$. Therefore, to choose the smoothing parameter $\lambda$, we simply propose to minimize the following empirical GCV-type criterion :
\begin{equation}
V(\lambda) = \frac{ \frac{1}{n} \sum_{i=1}^{n}(y_{i}-  \hat{f}_{n,\lambda}^{c}(x_{i}))^{2}} {  \int_{0}^{1} [Tr(I_{n} - A_{\lambda,t}) ]^{2} dt}. \label{eqGCV}
\end{equation}
In the above  equation, the quantity $ \int_{0}^{1} [Tr(I_{n} - A_{\lambda,t}) ]^{2} dt$ can be interpreted as a measure of the degree of freedom of the smoothing spline $ \hat{f}_{n,\lambda}^{c}$. The quantity $V(\lambda)$  is therefore the classical GCV criteria which is the ratio between the empirical error and the complexity of a smoothing procedure. To set a good penalization parameter $\lambda$, we simply use a grid search to minimize $V(\lambda)$.

\paragraph{Computational cost}
The proposed method has a relatively low computational cost compared to classical constrained optimization methods. If $t_{0},\ldots,t_{T-1}$ denotes a discretization of $[0,1]$ (with $T$ ndependent of $n$), our method requires for each $t_k$ the inversion of a symmetric definite matrix of size $n \times n$ which is possible using $\mathcal{O}(n^3)$ operations with a Cholesky algorithm for instance. The computational cost of our method is thus $\mathcal{O}(T n^3)$. Numerical computation of a constrained spline smoothing problem such as (\ref{eqsmoothsplineconst}) is generally done by using a set of inequality constraints to impose monotonicity on the value of the fitted function at a finite number of points. However, such algorithms  can be computationally intensive since solving a general problem of quadratic optimization with linear constraints is generally NP-hard (see e.g. Pardalos \& Vavasis \cite{nphard}). Primal-dual methods for instance can iteratively solve the problem but their complexity is larger than $\mathcal{O}(n^3)$. 


\section{A non-parametric regression problem under monotonicity constraints} \label{sec:stat}

Consider the standard nonparametric regression problem on a bounded
interval:
\begin{equation}
y_{i} = f(x_{i}) +  \epsilon_{i}, \; i=1,\ldots,n,  \label{eqreg}
\end{equation}
where $f : [0,1] \rightarrow \RR$  and $\epsilon_{i}$ are independent and identically distributed (i.i.d.) variables with zero mean and variance $\sigma^{2}$. The regression function $f$  is assumed to belong to a class of strictly increasing functions $\FF$ that satisfy some smoothness conditions to be defined later. Smoothing procedures for monotone regression can be found in   He \& Shi \cite{he}, Kelly \& Rice~\cite{kel}, Mammen~\cite{mam},
Mammen \& Thomas-Agnan~\cite{mam2}, Hall \&
Huang~\cite{hallhuang}, Mammen, Marron, Turlach \&
Wand~\cite{mam1}, Dette, Neumeyer \& Pilz~\cite{dette} and Antoniadis, Bigot \& Gijbels~\cite{antbig}.

In this section, we explain how homeomorphic splines can be used as a smooth and then monotonize method. Let  $\hat{f}_{n}$ be an unconstrained estimator obtained from the data $(y_{i},x_{i}), i =1,\ldots,n$ (e.g. by spline, kernel or wavelet smoothing). Our goal is to construct a monotone estimator $\hat{f}_{n}^{c}$ which inherits the asymptotic
properties of the unconstrained estimator $\hat{f}_{n}$ in terms of the
empirical mean squared error:
$
R_{n}(\hat{f}_{n}^{c},f) = \frac{1}{n}\sum_{i=1}^{n}(\hat{f}_{n}^{c}(x_{i})-f(x_{i}) )^{2}.
$
For this, starting from the values  $\hat{f}_{n}(x_i)$ instead of the observed $y_i$'s, take the vector field $v^{n,\lambda}$ which minimizes the following criterion:
$$
v^{n,\lambda} = \arg \min_{v \in \mathcal{X}} \int_{0}^{1} \frac{1}{n} \sum_{i=1}^{n}\left(\hat{f}^{n}(x_{i})-x_{i} -v_{t}(t\hat{f}^{n}(x_{i})+(1-t)x_{i}) \right)^{2}dt + \lambda \int_{0}^{1} \| h_{t} \|^{2}_{K} dt,
$$
where $v_{t}(x) = a_{1}^{t} + a_{2}^{t}x + h_{t}(x)$. Then $\hat{f}_{n,\lambda}^{c}$ is defined as the solution at time $t = 1$ of the ODE (\ref{eqodecont}) governed by the time-dependent vector field $v^{n,\lambda}$.

\subsection{Asymptotic properties of the monotone estimator $\hat{f}_{n,\lambda}^{c}$}

The following theorem shows that under mild conditions on the design and the unconstrained estimator, the monotone estimator $\hat{f}_{n,\lambda}^{c}$ inherits the asymptotic properties of $\hat{f}^{n}$ in term of rate of convergence. To the best of our knowledge, this is the first consistency result on estimators defined through large diffeomorphism models governed by ODE.\\

\begin{theo} \label{theomain}
Assume that the conditions of Theorem \ref{thode} are satisfied, and that the kernel $K: \RR^{2} \to \RR$ is continuous. Moreover assume that the function $f$ is continuously differentiable on $[0,1]$ with $f'(x) > 0$ for
all $x \in [0,1]$ and that there exists a time-dependent vector field $v^{f} \in
\mathcal{X}^{2}$ such that for all $t \in [0,1]$:
$$
v^{f}_{t}( \phi_{t}(x) ) = f(x) - x \mbox{ for all } x \in [0,1],
$$
where $\phi_{t}(x) = tf(x) + (1-t)x$. Suppose that the unconstrained estimator  $\hat{f}^{n}$ and the points $x_{i},i=1,\ldots,n$ satisfy the following property:\\

\begin{description}
\item[A1] for all $t \in [0,1]$ and all $1 \leq i,j \leq n$ with $i \neq j$,
\begin{equation}
t\hat{f}^{n}(x_{i})+(1-t)x_{i} \neq t\hat{f}^{n}(x_{j})+(1-t)x_{j} \; a.s. \label{eqalmost}
\end{equation}
\end{description}

\noindent Then, $v^{n,\lambda} \in \mathcal{X}^{2} a.s. $ and thus $\hat{f}_{n,\lambda}^{c}$ is a monotone function on $[0,1]$. If we further assume that:\\

\begin{description}
\item[A2] there exists a weight function $\omega : [0,1] \to ]0,+\infty[$ such that for any  $g \in \mathcal{C}^1([0,1],\RR) $ one has $\lim_{n \to +\infty} \frac{1}{n} \sum_{i=1}^{n} g(x_{i}) \to \int_{0}^{1} g(x) \omega(x) dx$,
\item[A3] $R_{n}(\hat{f}_{n},f) \to 0 \mbox{ in probability as } n \to +\infty$.\\
\end{description}

\noindent Then, for any sequence $\lambda = \lambda_{n} \to 0$, we have that there exists a deterministic constant $\Lambda_{1}$ (not depending on $n$)   such that with probability tending to one as $n \to +\infty$:
$$
R_{n}(\hat{f}_{n,\lambda_{n}}^{c},f) \leq \Lambda_{1} \left( R_{n}(\hat{f}_{n},f) +  \lambda_{n} \right).
$$
\end{theo}

Equation (\ref{eqalmost}) may not be satisfied for time points $t_{ij}$ such that $\frac{1-t_{ij}}{t_{ij}} = \frac{\hat{f}^{n}(x_{i})- \hat{f}^{n}(x_{j})}{x_{j}-x_{i}}$. Since the function $t \to \frac{1-t}{t}$ is injective, this can only happen for a finite number of time points $t \in [0,1]$. Hence assumption  {\bf A1} is generally satisfied provided the design points are distinct. Moreover, if equation (\ref{eqalmost}) is not satisfied for some points $t$, one can argue that it is possible to modify the estimator $\hat{f}^{n}$ without changing its asymptotic properties (by slightly varying e.g. the smoothing parameter used to compute it) such that (\ref{eqalmost}) is true for any $t \in [0,1]$ and all $1 \leq i,j \leq n$. Note that  under  assumption  {\bf A1},  Proposition \ref{propexist} implies that $v^{n,\lambda}$ can be easily implemented using unconstrained spline smoothing.

The assumption {\bf A2} means that in some sense the design points are sampled according to the density $\omega(x)$. The assumption {\bf A3} is satisfied whenever the expected empirical mean squared error $\EE R_{n}(\hat{f}_{n},f)$ converges to zero as $n \to +\infty$. The fact that $v^{n,\lambda} \in \mathcal{X}^{2}$ guarantees that $\hat{f}_{n}$ is a monotone function. Moreover, one can see that if  $\lambda_{n}$ decays as fast as the empirical error  $R_{n}(\hat{f}_{n},f)$ then the estimator $\hat{f}_{n,\lambda_{n}}^{c}$ has the same asymptotic convergence rate than the unconstrained estimator $\hat{f}_{n}$.
Similar results for smooth and then monotonize approaches are discussed in Mammen \& Thomas-Agnan~\cite{mam2},  Mammen, Marron, Turlach \& Wand~\cite{mam1}. However, the advantages of our approach over existing methods are the following: it yields a monotone smoother which has a closed form expression and which is guaranteed to be monotone on the whole interval $[0,1]$. Moreover, our approach is very general as it is not restricted to functions belonging to Sobolev spaces. 

\subsection{Optimal rate of convergence for Sobolev spaces}

Let us return to the specific case where $f \in H^{m}([0,1])$ and $\HH_{K} = \HH^{m}(\RR)$. The asymptotic properties and optimal rates of convergence (in the minimax sense) of unconstrained estimators for functions belonging to Sobolev spaces has been extensively studied (see e.g. Nussbaum \cite{nuss}, Speckman \cite{spec}).  The estimator $\hat{f}_{n,\gamma}$ of Speckman \cite{spec} is based on the use of the Demmler-Reinsch spline basis and on a smoothing parameter $\gamma$ (see Speckman \cite{spec} for further details). Speckman \cite{spec}  has shown that for an appropriate choice $\gamma^{\ast}$ then $ \EE R_{n}(\hat{f}_{n,\gamma^{\ast}},f) =\opO \left( n^{-\frac{2m}{2m+1}}\right)$ if $f \in H^{m}([0,1])$ which is known to be the minimax rate of convergence for functions belonging to Sobolev balls. This result is based on the assumption that the design points  are such that $x_{i} = G((2i-1)/2n)$ where $G : [0,1] \to [0,1]$ is a continuously differentiable function with $G'(x) \geq c > 0$ for some constant $c$.  Hence, the estimator of Speckman \cite{spec}  satisfies Assumption  {\bf A2} with $\omega(x) = \frac{1}{G'(G^{-1}(x))}$, and one can check that Assumption  {\bf A1} also holds. The following corollary is thus an immediate consequence of Theorem \ref{theovf} and Theorem \ref{theomain}:\\

\begin{coro}
Assume that $\tilde{\HH} = Span\{1,x\} + \HH^{m}(\RR)$. Let $m \geq 2$ and $f  \in   H^{m}([0,1])$ be such that $f'(x) > 0$ for all $x \in [0,1]$. Then, the monotone estimator $\hat{f}_{n,\lambda_{n}}^{c}$ based on the minimax estimator $\hat{f}_{n,\gamma^{\ast}}$ of Speckman \cite{spec}  is such that $ R_{n}(\hat{f}_{n,\lambda_{n}}^{c},f) =\opO_{P}\left( n^{-\frac{2m}{2m+1}}\right),$ provided $\lambda_{n} = \opO(n^{-\frac{2m}{2m+1}})$.\\
\end{coro}

To obtain an adaptive choice of $\lambda$ (not depending on the unknown regularity $m$ of $f$), the above Corollary suggests to take $\lambda_{n} = \frac{1}{n}$ to have a monotone estimator whose empirical mean squared error decays as fast as $R_{n}(\hat{f}_{n},f)$. This choice may yield satisfactory estimates but in our simulations a data-based choice  for $\lambda$ using a GCV criteria (\ref{eqGCV})  gives much better results.


\section{Numerical experiments} \label{secsimu}

\subsection{1D case and monotonicity}

Dette, Neumeyer and Pilz~\cite{dette} have recently proposed another type of smooth and then monotonize method which combines density and regression estimation with kernel smoothers. This approach has been shown to be very successful on many simulated and real data sets (see Dette and Pilz~\cite{detpilz}) and we shall therefore use it as a benchmark to assess the quality of our monotone estimator.  Similarly to our approach, it requires a preliminary unconstrained  estimator  $\hat{f}_{n}$. This estimator is then used to estimate the inverse $f^{-1}$ of the regression function. For this Dette, Neumeyer and Pilz~\cite{dette} propose to use the following estimator
$$
\hat{m}^{-1}_{n}(x) = \frac{1}{N h_{d}} \sum_{i=1}^{N} \int_{-\infty}^{t} K_{d} \left( \frac{\hat{f}_{n}(\frac{i}{N})-u}{h_{d}} \right)du,
$$
where $K_{d}$ is a positive kernel function with compact support, $h_{d}$ a bandwidth that controls the smoothness of $\hat{m}^{-1}_{n}$ and $N$ is an integer not necessarily equal to the sample size $n$ which controls the numerical precision of the procedure. A monotone estimator $\hat{m}_{n}$ is then obtained by reflection of the function $\hat{m}^{-1}_{n}$ at the line $y=x$ (see Dette, Neumeyer and Pilz~\cite{dette} for further details). In Dette and Pilz~\cite{dette}, it is proposed to use a local linear estimate (see Fan and Gijbels \cite{fangij}) with Epanechnikov kernel for the unconstrained estimator $\hat{f}_{n}$. The bandwidth $h_{r}$ of this unconstrained estimate is chosen as $\hat{h}_{r} = \left(\frac{\hat{\sigma}^{2}}{n} \right)^{1/5}$, where $\hat{\sigma}^{2}  =  \frac{1}{2(n-1)} \sum_{i=1}^{n-1}\left(y_{(i+1)}-y_{(i)} \right)^{2}$.  For the choice of the bandwidth $h_{d}$, it is recommended to choose $h_{d} = h_{r}^{3}$. However, for a fair comparison with our data-based choice of $\lambda$ via GCV, the best choice for $h_{d}$ is chosen by cross-validation via a grid search.

\begin{figure}[htdp]
\centering \subfigure[] { \includegraphics[width=3cm]{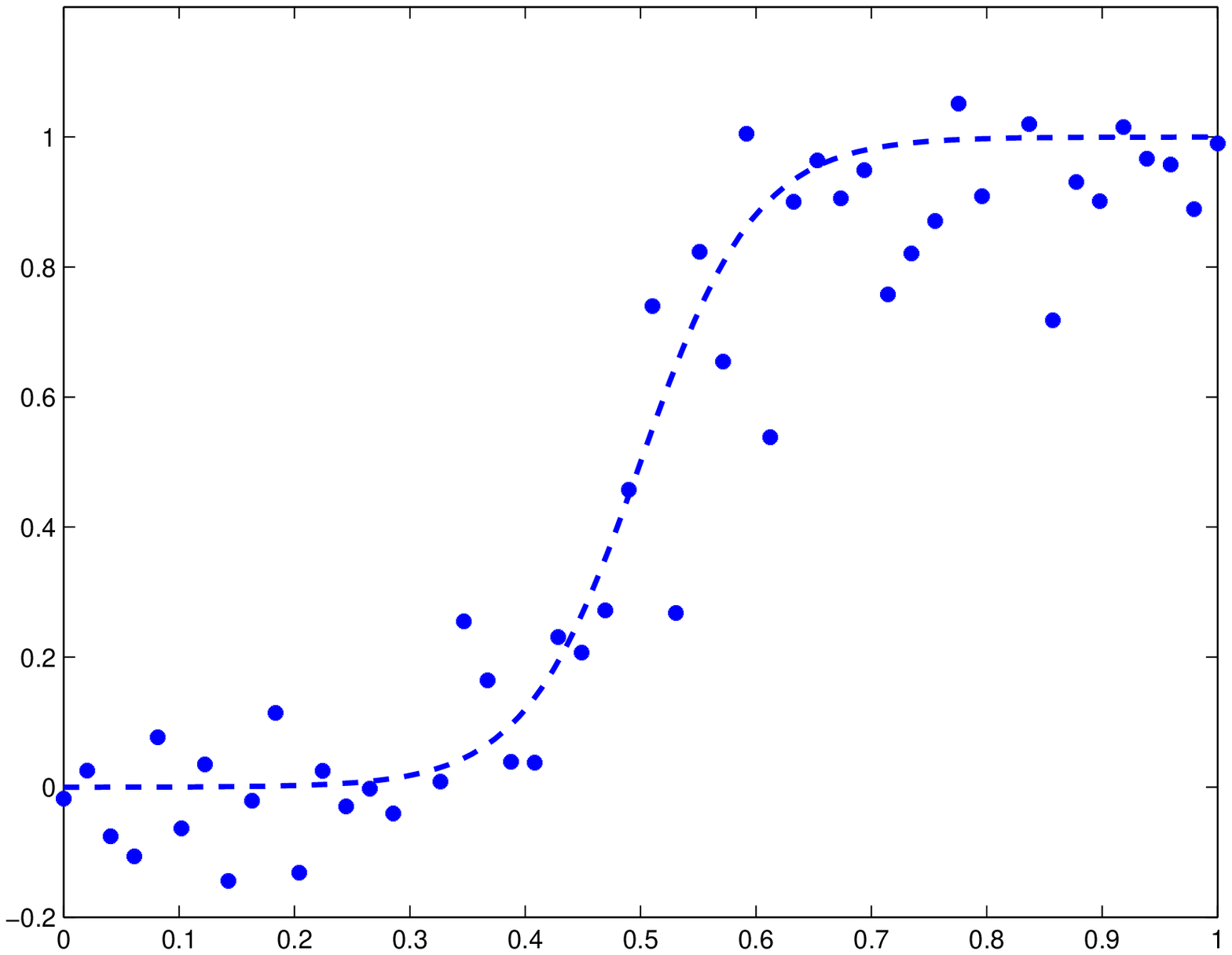} }
\subfigure[] { \includegraphics[width=3cm]{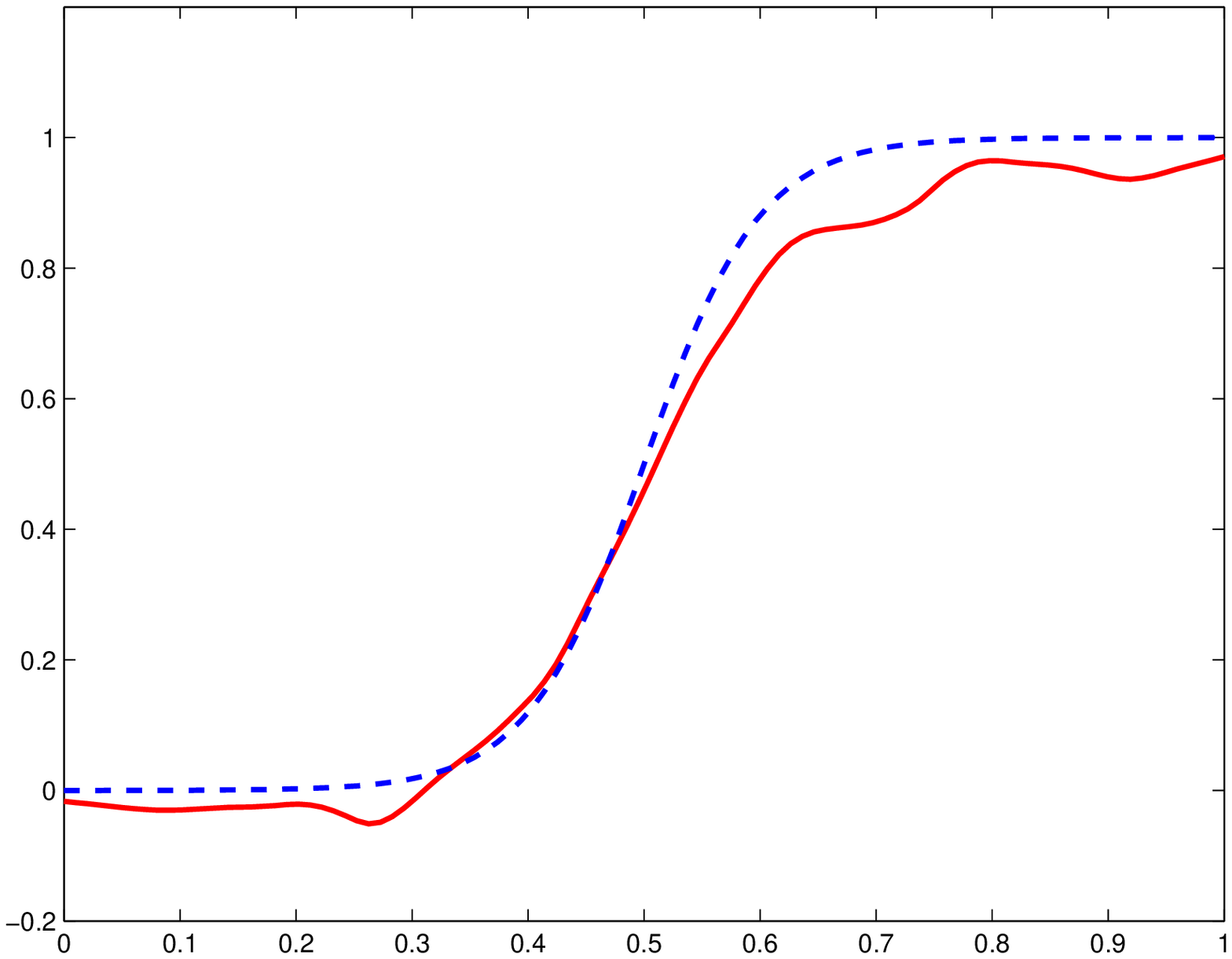} }
\subfigure[] { \includegraphics[width=3cm]{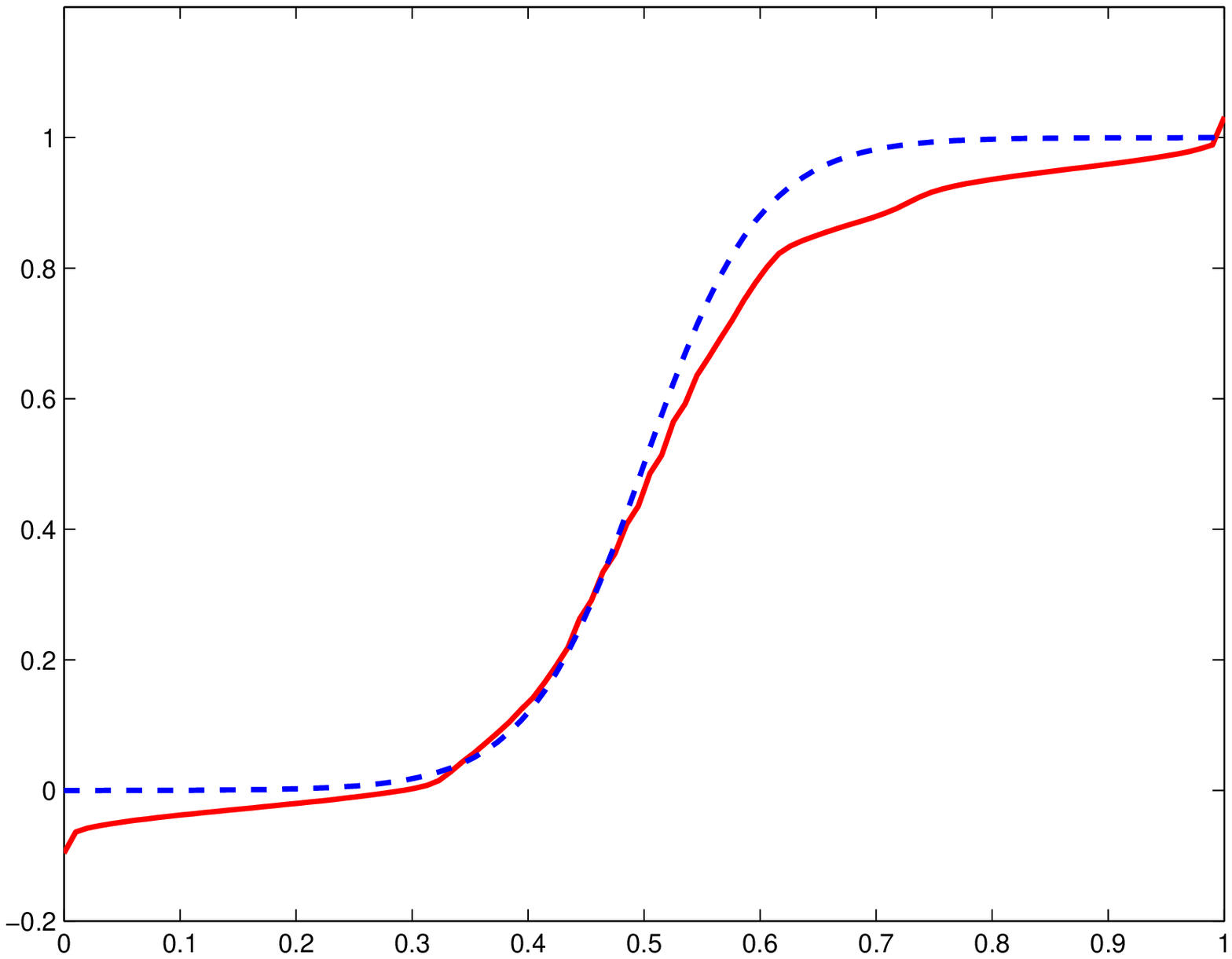} }
\subfigure[] {\includegraphics[width=3cm]{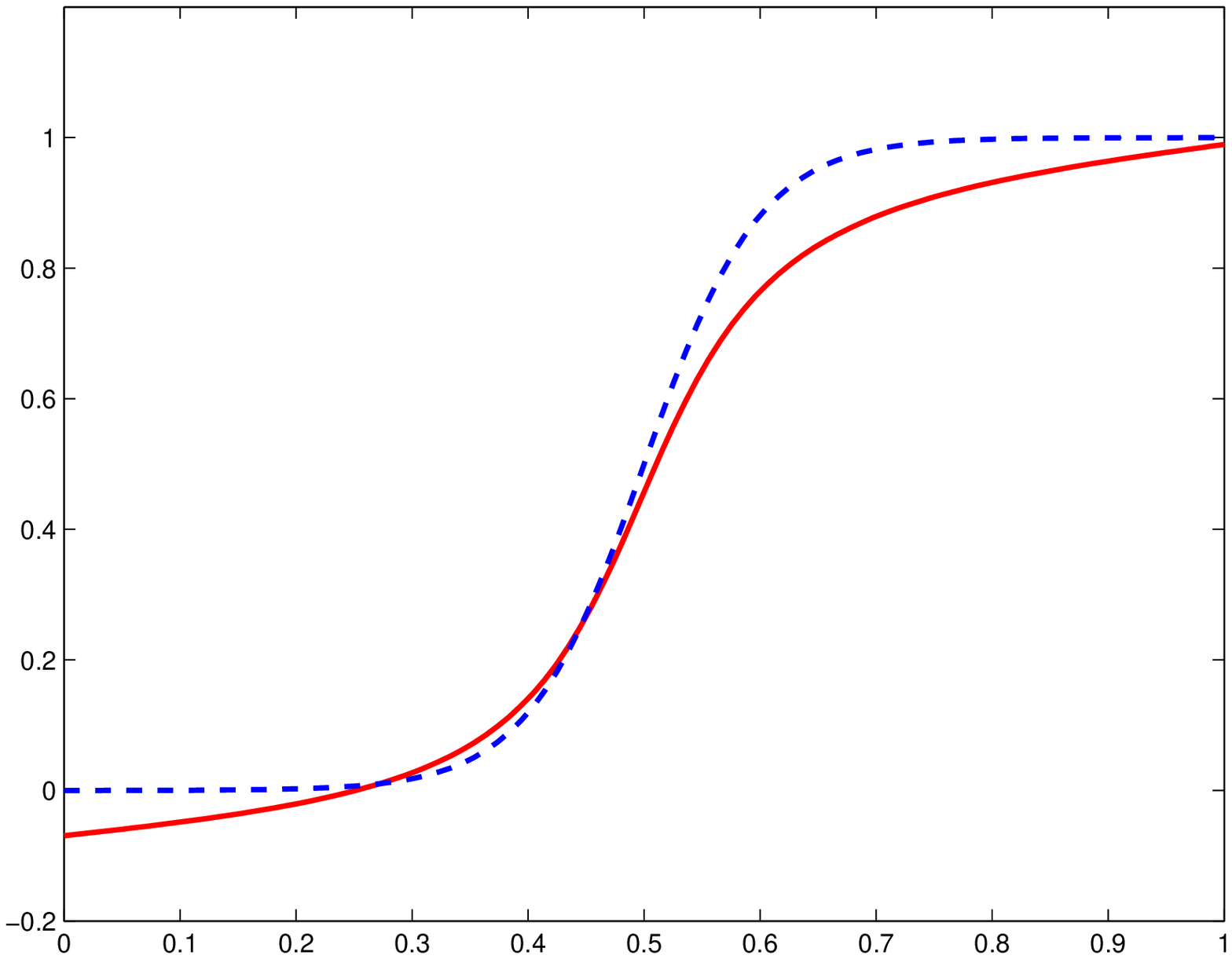} }

\caption{{\em {\small Signal $m_{1}$: the dotted line is the unknown regression
function, (a) noisy data  with $SNR=3$, (b) Local Linear Unconstrained
Estimator, (c) Dette {et al.}'s estimator, (d) Homeomorphic
Smoothing Spline based on the Local Linear estimator.}}}
\label{figm1}
\end{figure}

\begin{figure}[htdp]
\centering \subfigure[] { \includegraphics[width=3cm]{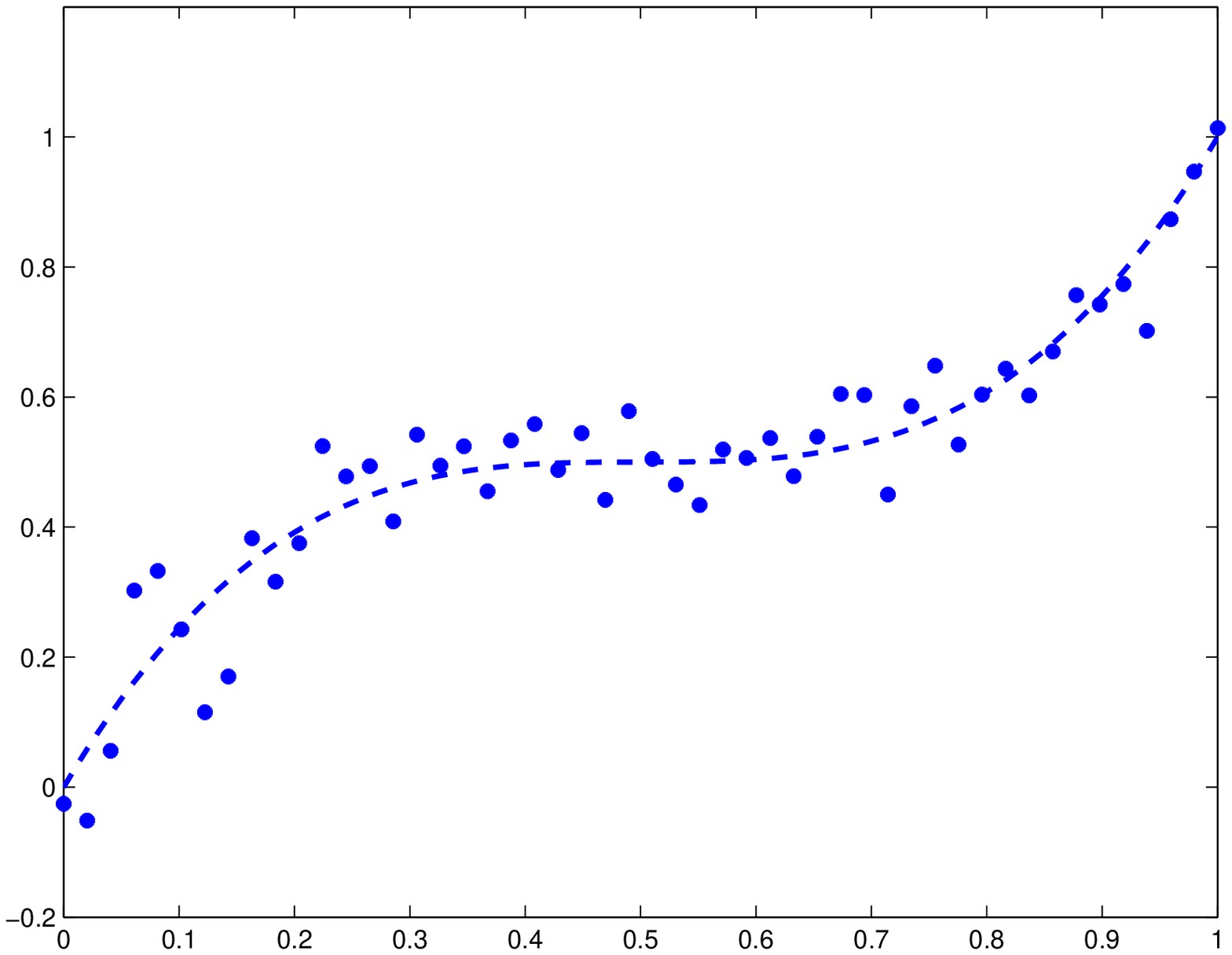} }
\subfigure[] { \includegraphics[width=3cm]{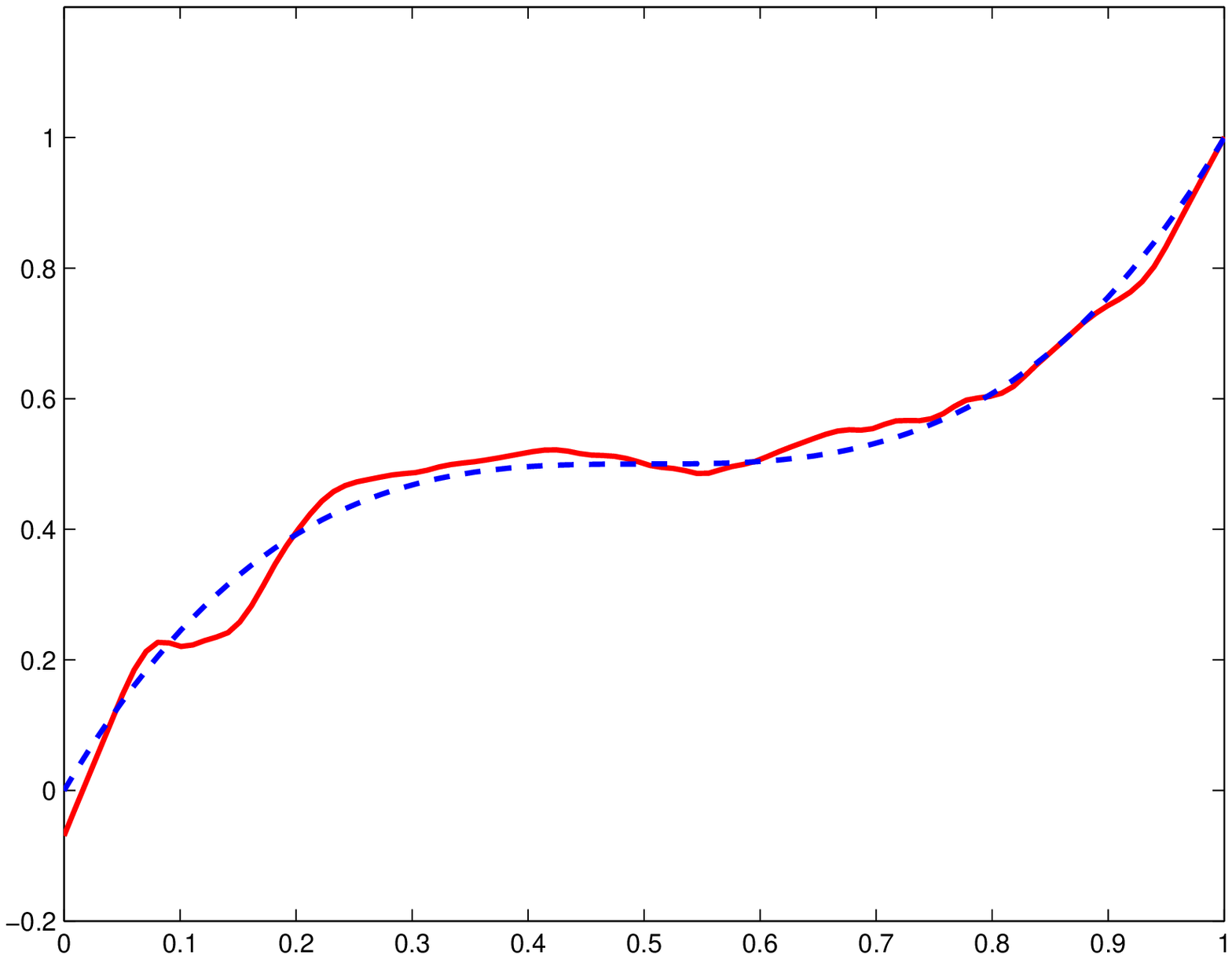} }
\subfigure[] { \includegraphics[width=3cm]{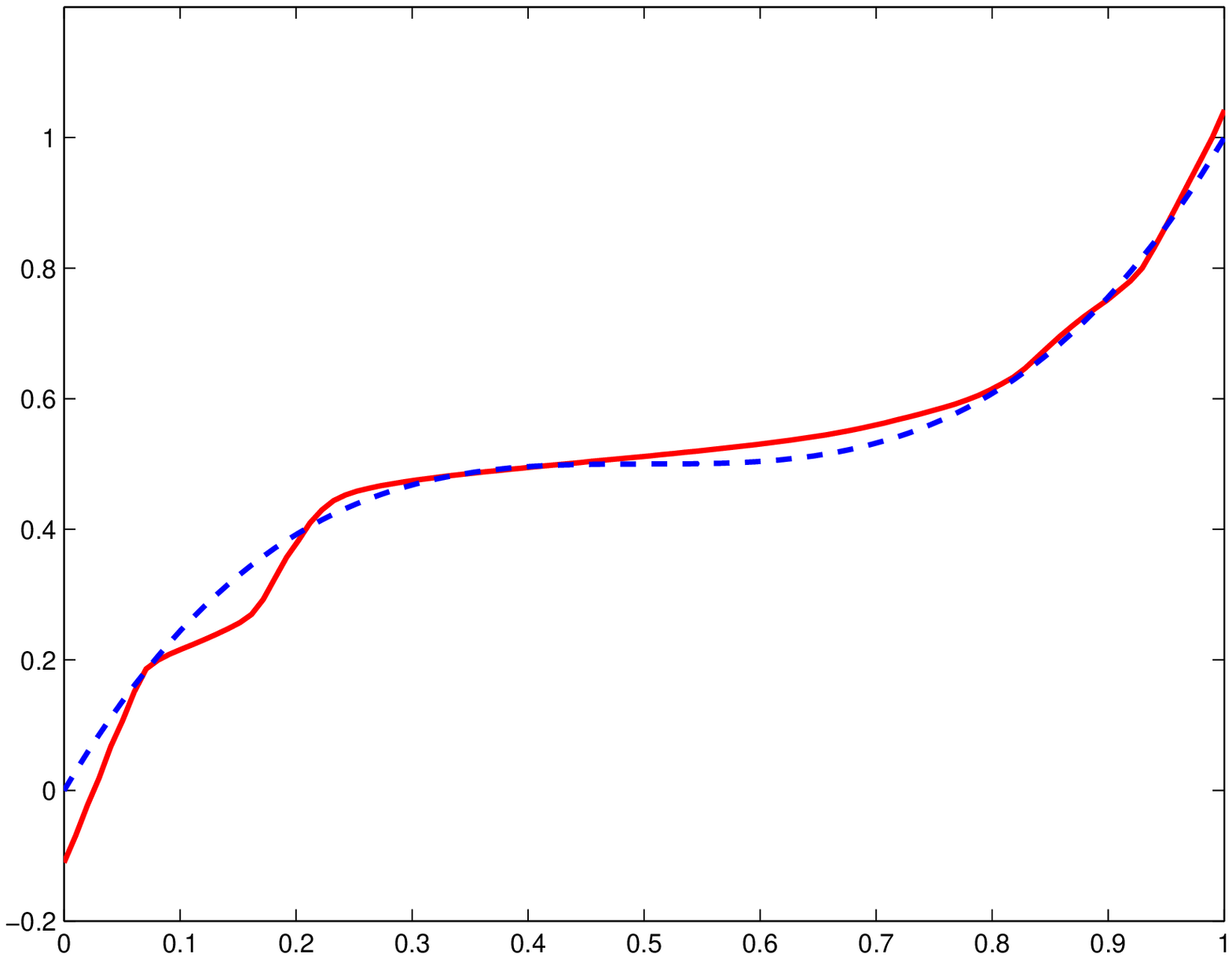} }
\subfigure[] {\includegraphics[width=3cm]{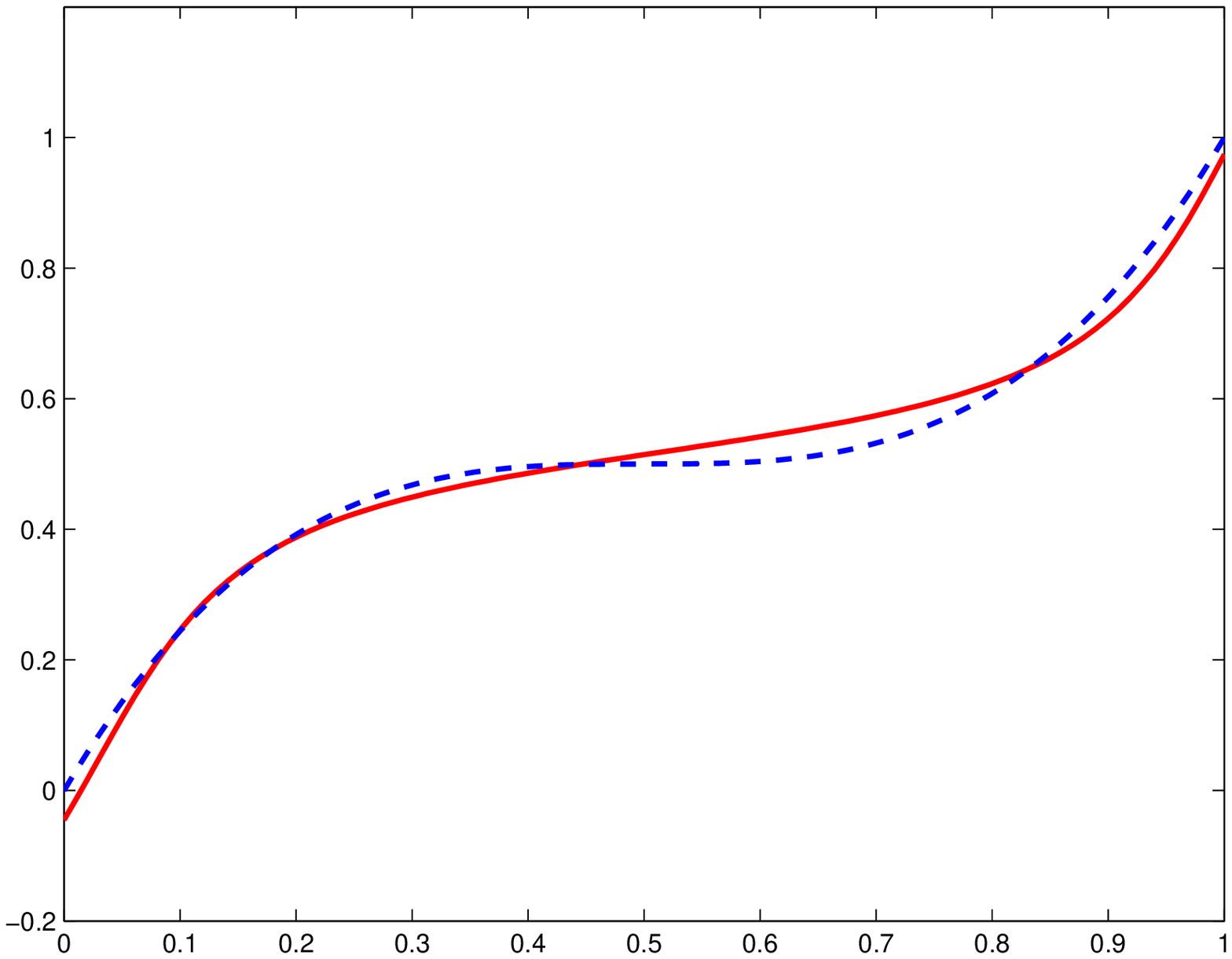} }

\caption{{\em {\small Signal $m_{2}$: the dotted line is the unknown regression
function, (a) noisy data  with $SNR=3$, (b) Local Linear Unconstrained
Estimator, (c) Dette {et al.}'s estimator, (d) Homeomorphic
Smoothing Spline based on the Local Linear estimator.}}}
\label{figm2}
\end{figure}

\begin{figure}[htdp]
\centering \subfigure[] { \includegraphics[width=3cm]{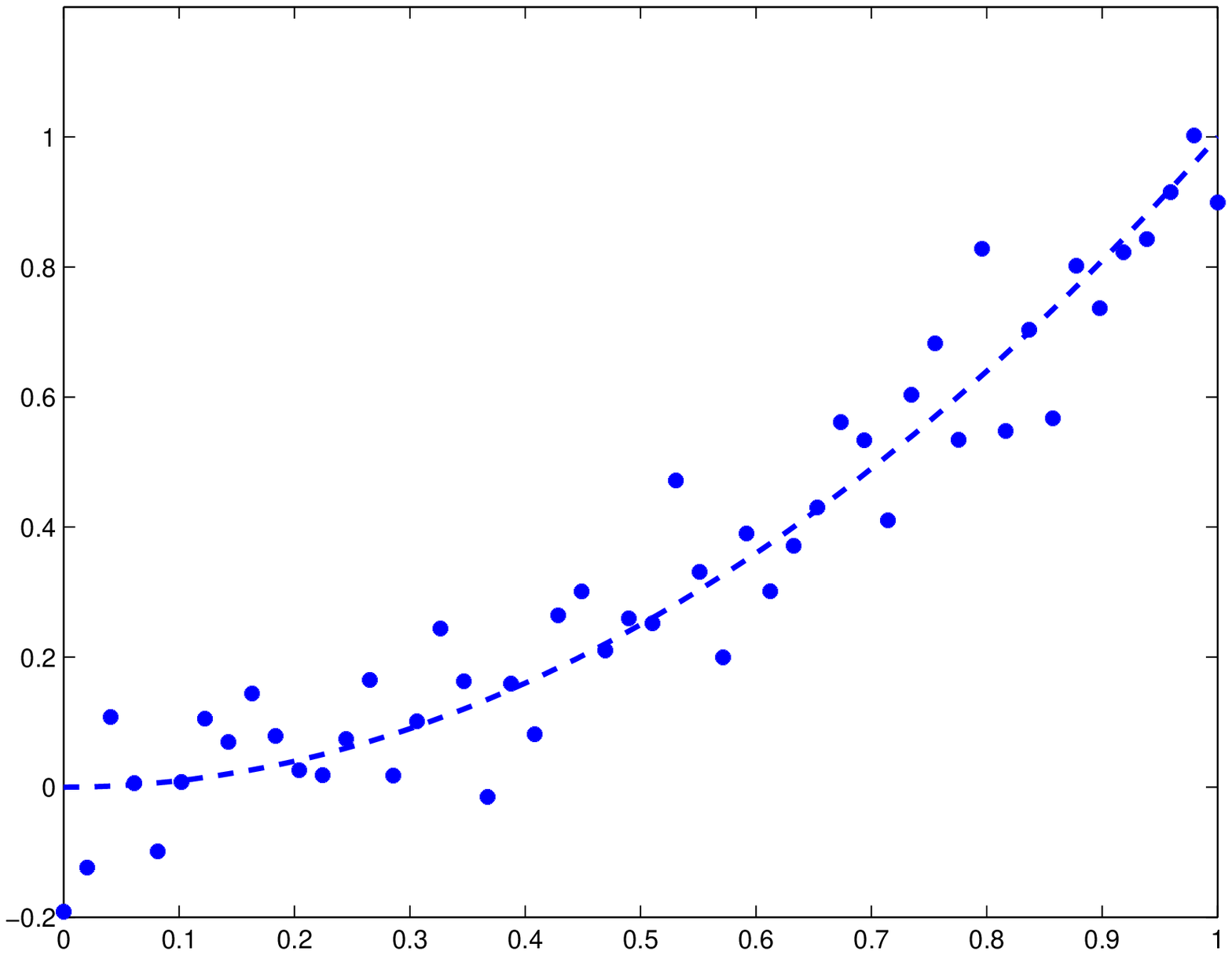} }
\subfigure[] { \includegraphics[width=3cm]{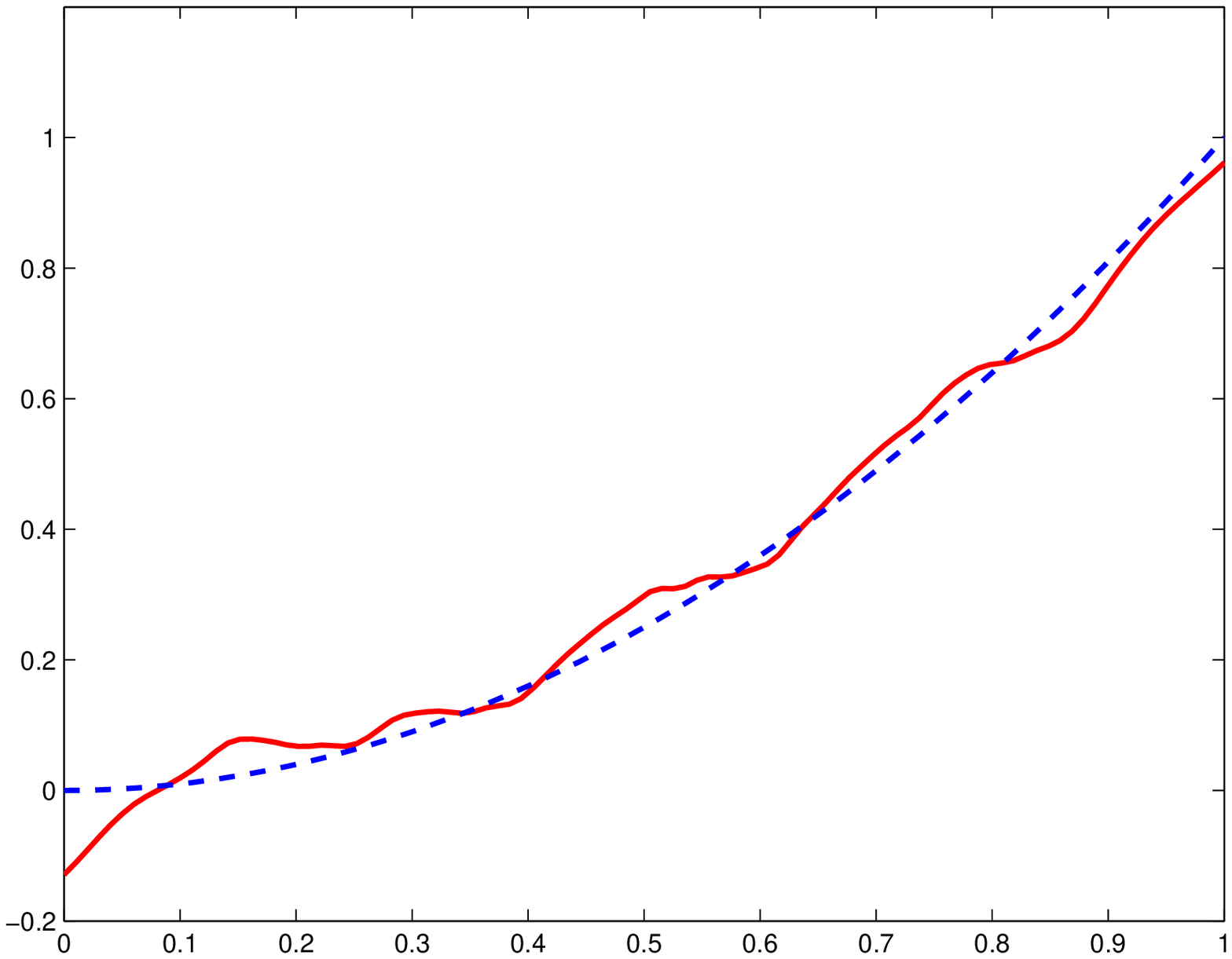} }
\subfigure[] { \includegraphics[width=3cm]{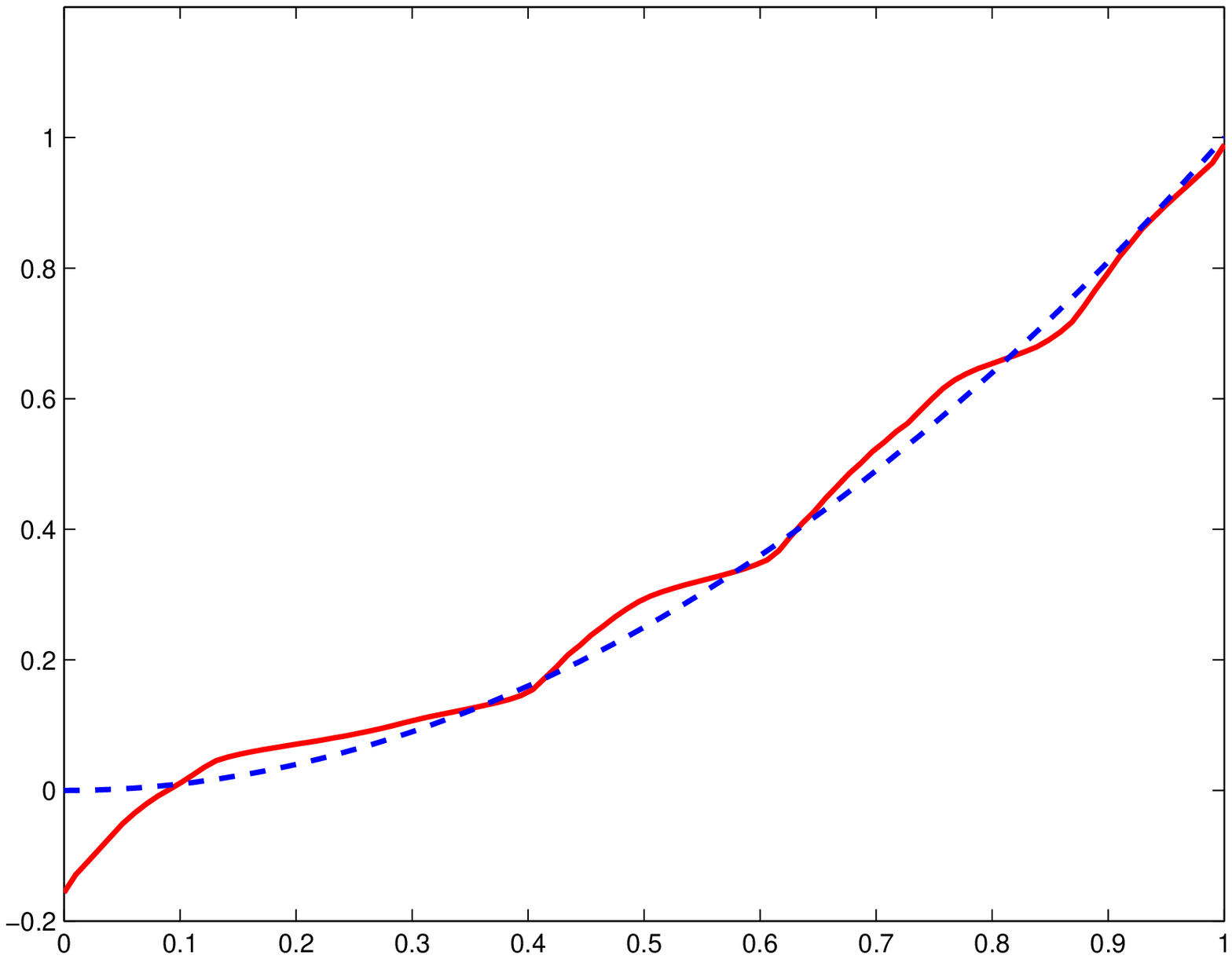} }
\subfigure[] {\includegraphics[width=3cm]{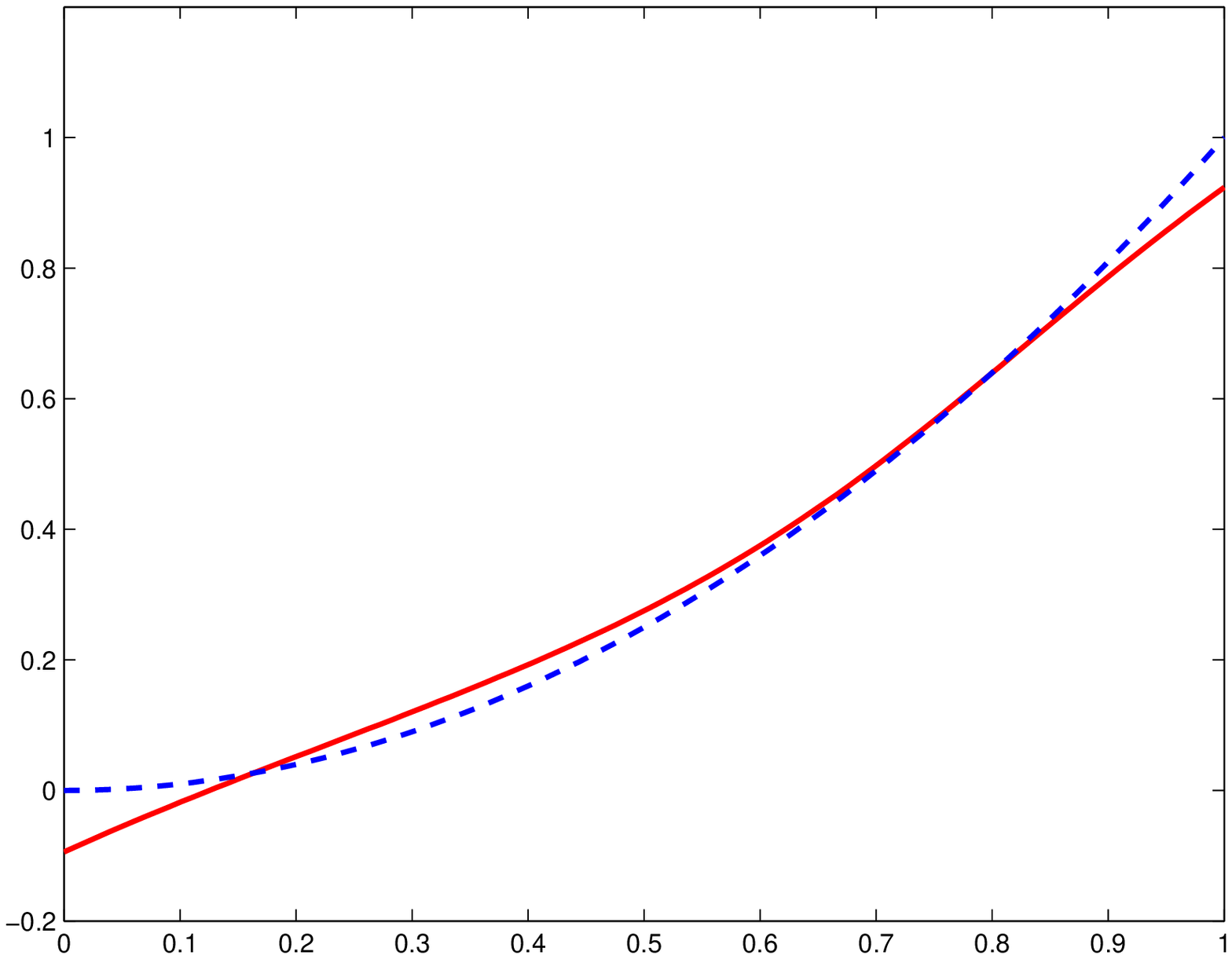} }

\caption{{\em {\small Signal $m_{3}$: the dotted line is the unknown regression
function, (a) noisy data  with $SNR=3$, (b) Local Linear Unconstrained
Estimator, (c) Dette {et al.}'s estimator, (d) Homeomorphic
Smoothing Spline based on the Local Linear estimator.}}}
\label{figm3}
\end{figure}

We investigate the regression model with a regular design i.e. $x_{i} = \frac{i}{n}, i=1,\ldots,n$,  normally distributed errors, sample size $n = 50$ and a signal to noise ratio (SNR) of 3.  The signal-to-noise ratio is measured as $\hbox{sd}(f(x))/\sigma$,
where $\hbox{sd}(f(x))$ is the estimated standard deviation of the regression function, $f(x_{i}) $ over the sample $i = 1, \dots, n$, and $\sigma$  is the true standard deviation of the noise in the data. The monotone regression functions that we consider are (see  Dette and Pilz~\cite{dette})
$$
m_{1}(x) = \frac{\exp(20(x-1/2))}{1 + \exp(20(x-1/2))}, \; m_{2}(x) = \frac{1}{2}(2x-1)^{3}+\frac{1}{2}, \; m_{3}(x) = x^{2}.
$$
These functions  correspond to, respectively,  a function with a  ``continuous jump'', a strictly increasing curve with a plateau, and a convex function. The different functions are
displayed in Figures \ref{figm1}-\ref{figm3}.

In Figures \ref{figm1}-\ref{figm3}, we present some curves for the estimates of these three test functions. A Gaussian kernel has been used to compute the homeomorphic smoothing splines (using other kernels gives similar results). For the choice of the regularization parameter $\lambda$ the GCV criterion (\ref{eqGCV}) is used, and recall that we use cross-validation for the choice of $h_{d}$. As one can see in Figures \ref{figm1}-\ref{figm3}, the homeomorphic smoothing spline based on the  local linear estimator gives results similar to those obtained via the estimator of Dette, Neumeyer and Pilz~\cite{dette}. However, our approach yields  monotone estimator that are visually much smoother and very close to the true regression in all cases. Homeomorphic smoothing splines also seems to give very nice results even if the unconstrained estimator $\hat{f}_{n}$ is very oscillating as it is the case for the local linear estimator in Figure \ref{figm2} and  Figure \ref{figm3}.


\begin{figure}[htdp]
\centering \subfigure[] { \includegraphics[width=3cm]{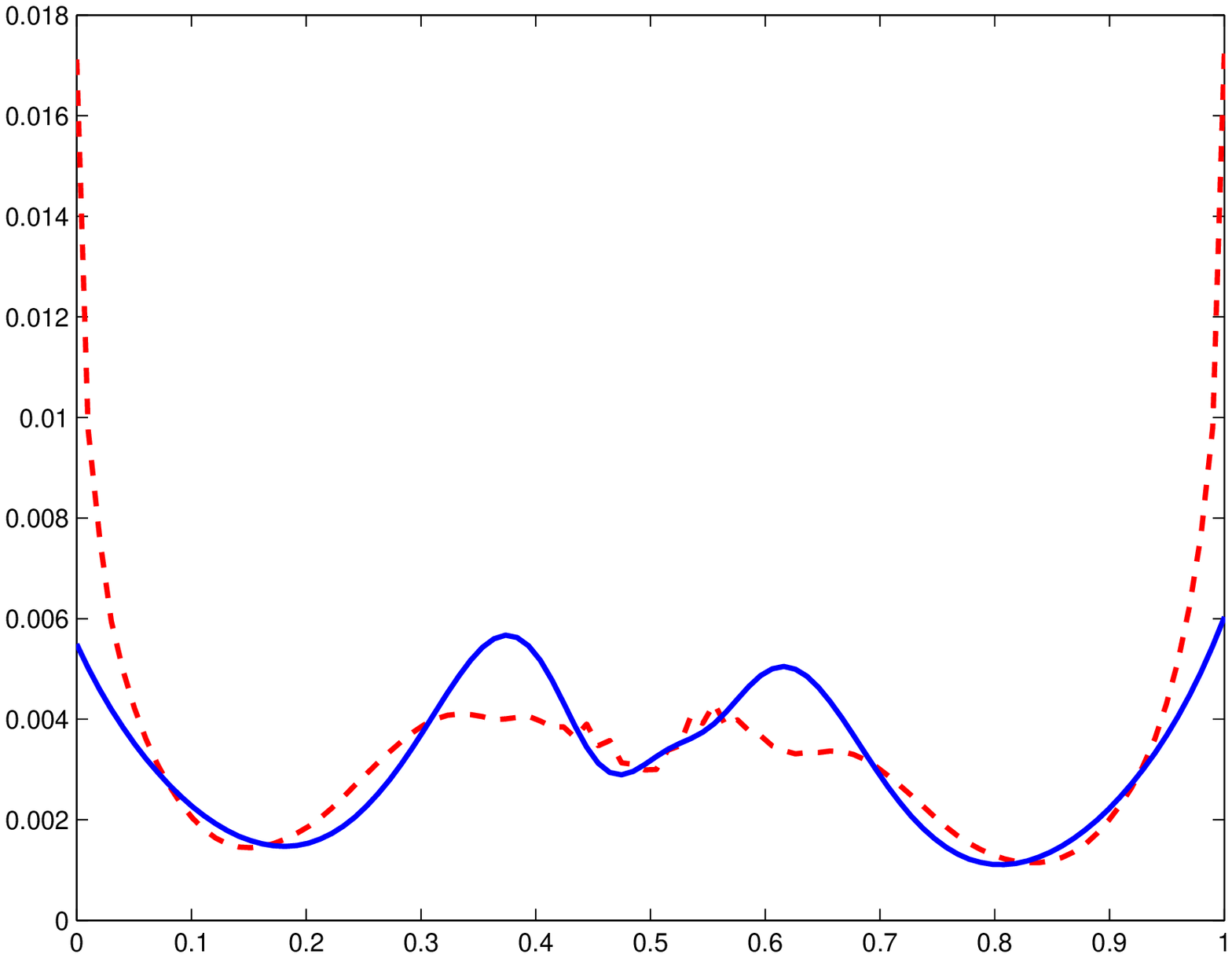} }
\subfigure[] { \includegraphics[width=3cm]{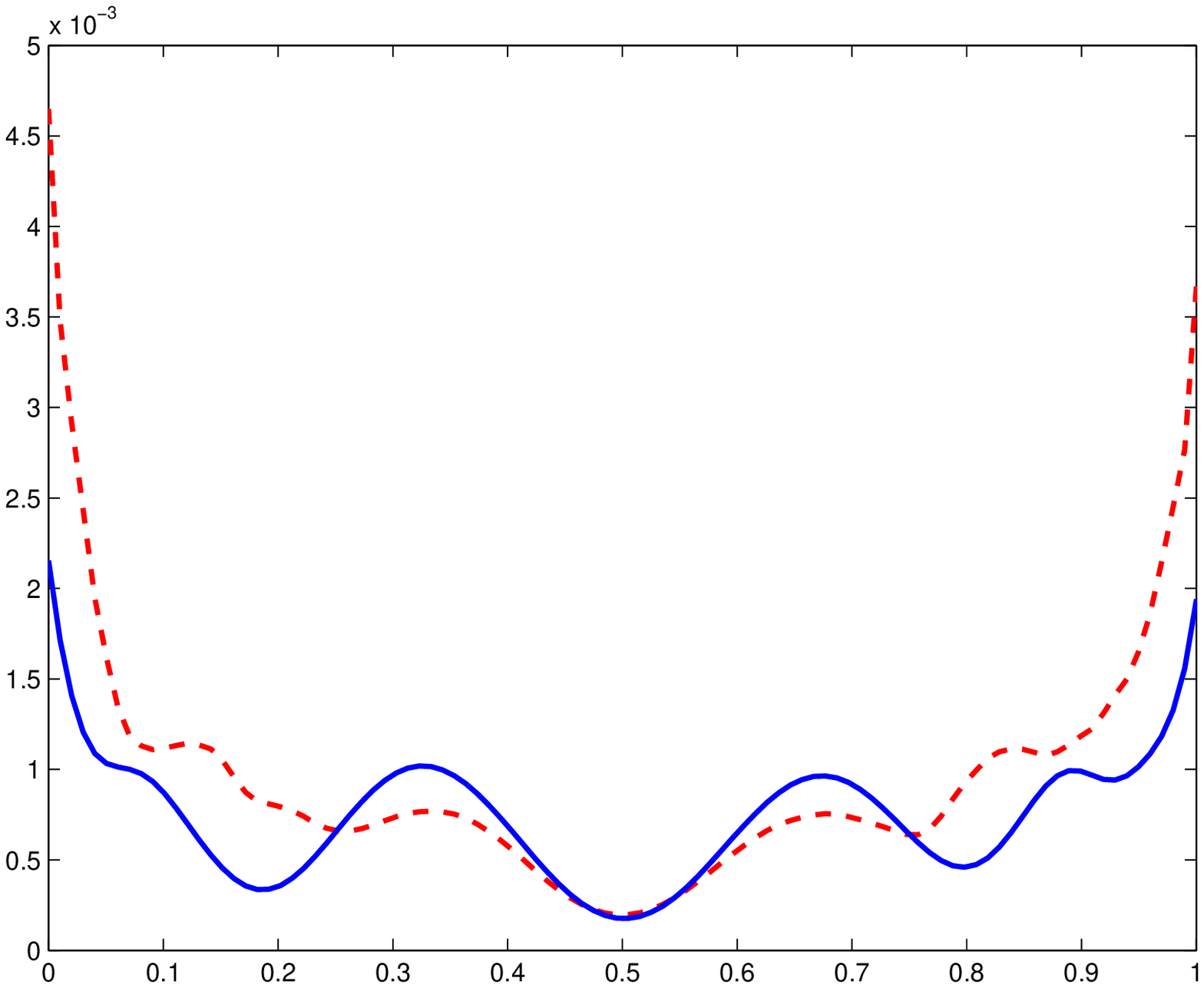} }
\subfigure[] { \includegraphics[width=3cm]{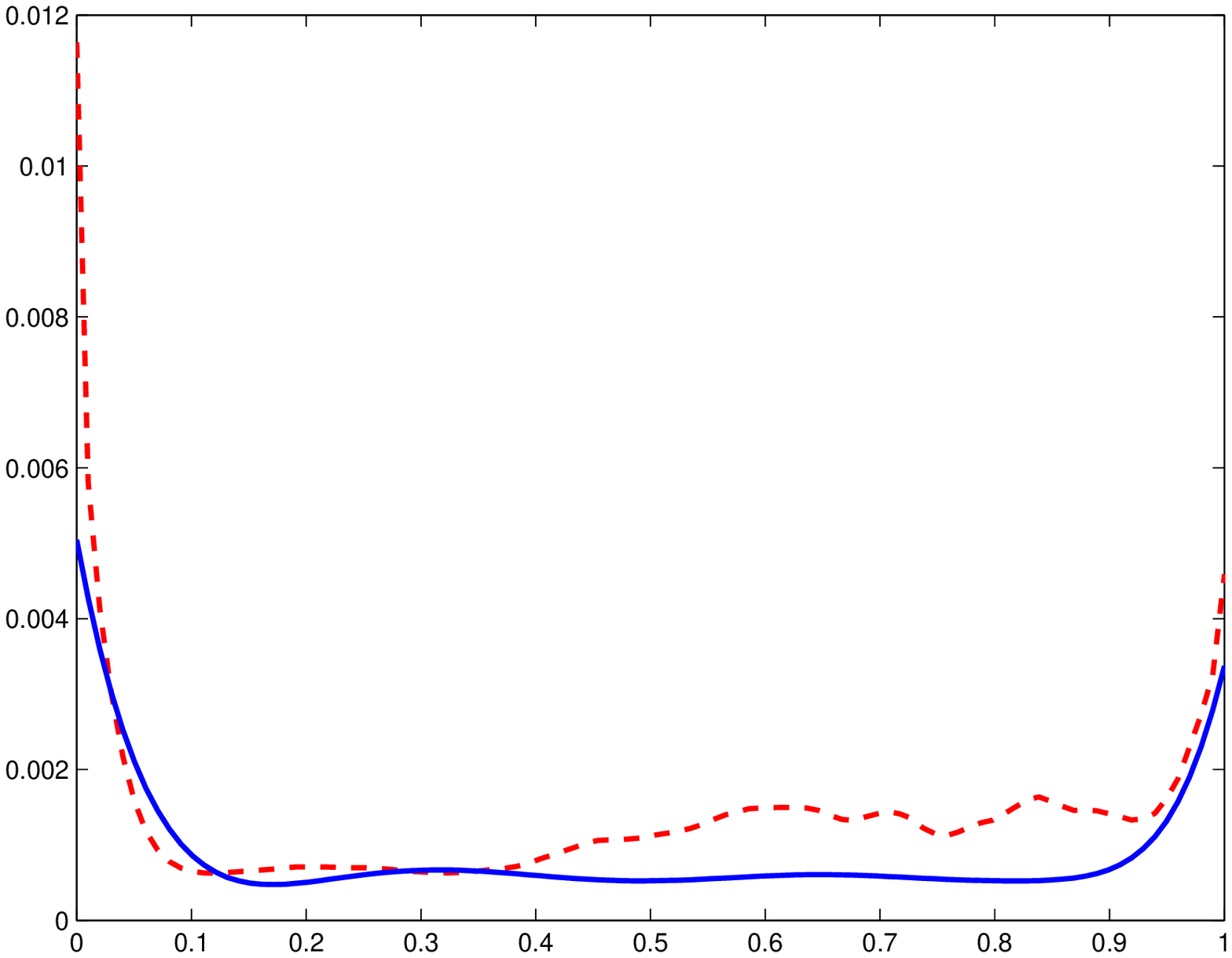} }

\caption{{\em  {\small Simulated mean squared error    with $SNR=3$ (computed
over 100 simulations runs) on signal $m_1,m_2$ and $m_3$: Dette {\it et al.}'s estimator (dashed
curves) and homeomorphic smoothing spline (solid curves)  for the three regression functions.}}} \label{figMSE3}
\end{figure}

\begin{table}[htdp]
\caption{Mean integrated squared error (MISE) over the 100 simulations for each method.}
\begin{small}
\begin{center}
\begin{tabular}{c|c|c|c}
& Signal $m_{1}$ & Signal $m_{2}$ & Signal $m_{3}$   \\ \hline

Homeomorphic smoothing spline &   0.0032 & 0.00076 & 0.00089 \\  \hline

Dette {\it et al.}'s estimator  &   0.0035 & 0.00098 & 0.0014 \\

\end{tabular}
\end{center}
\label{tabMISE}
\end{small}
\end{table}


To compare  these two monotone estimates, we have used 100 simulations runs for each regression function. The same unconstrained estimator (a local linear estimate with Epanechnikov kernel) is used. For the 100 simulations, we have calculated the pointwise mean squared error (MSE) for the two estimates $\hat{f}_{n}^{c}$ and $\hat{m}_{n}$, evaluated on an equidistant grid of size $2n$. Curves for the MSE of the three estimate are displayed in Figure \ref{figMSE3}. Again, these simulations clearly show that our approach compares similarly to the monotone estimator of Dette, Neumeyer \& Pilz~\cite{dette} for the signal $m_{1}$ and $m_{2}$, and outperforms Dette {et al.}'s estimator for the function $m_{3}$. Table \ref{tabMISE} shows that it gives better results in terms of mean integrated  squared error (MISE over $[0,1]$) for the three test functions.

\subsection{2D experiments and diffeormorphic matching} \label{sec:H2}
\begin{figure}[htdp]
\centering \subfigure[] { \includegraphics[width=3cm]{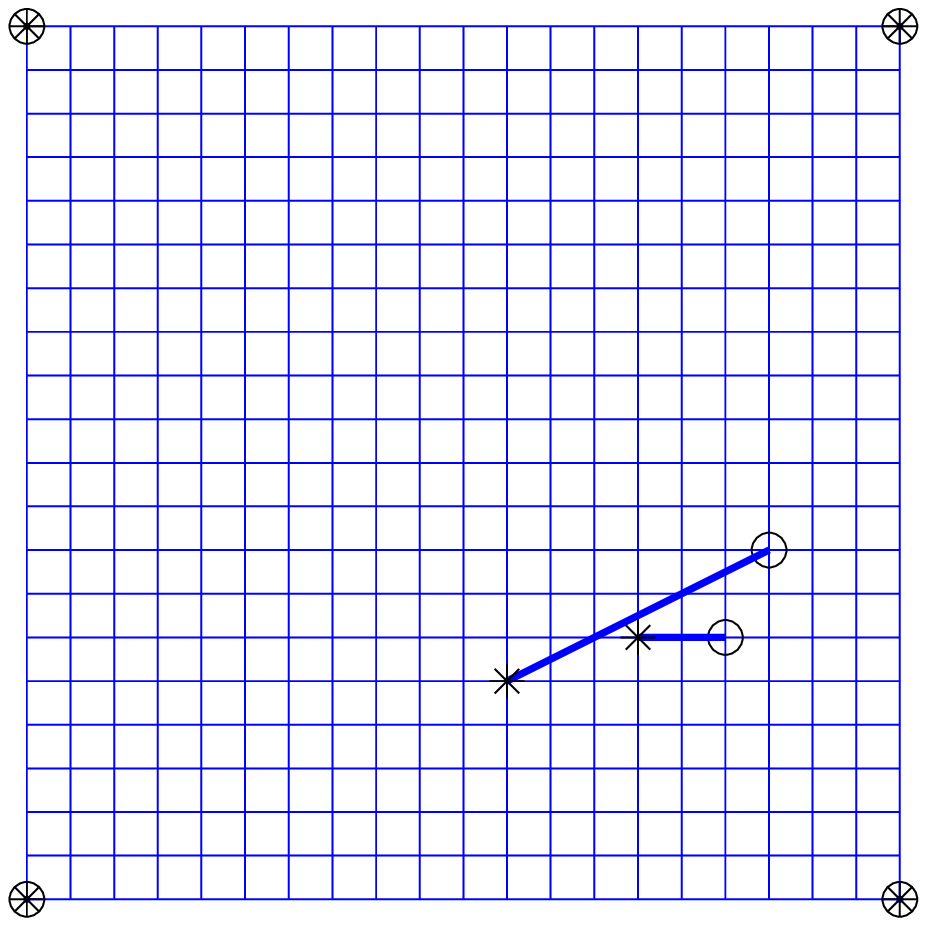} } \hspace{0.5cm}
\subfigure[] { \includegraphics[width=3cm]{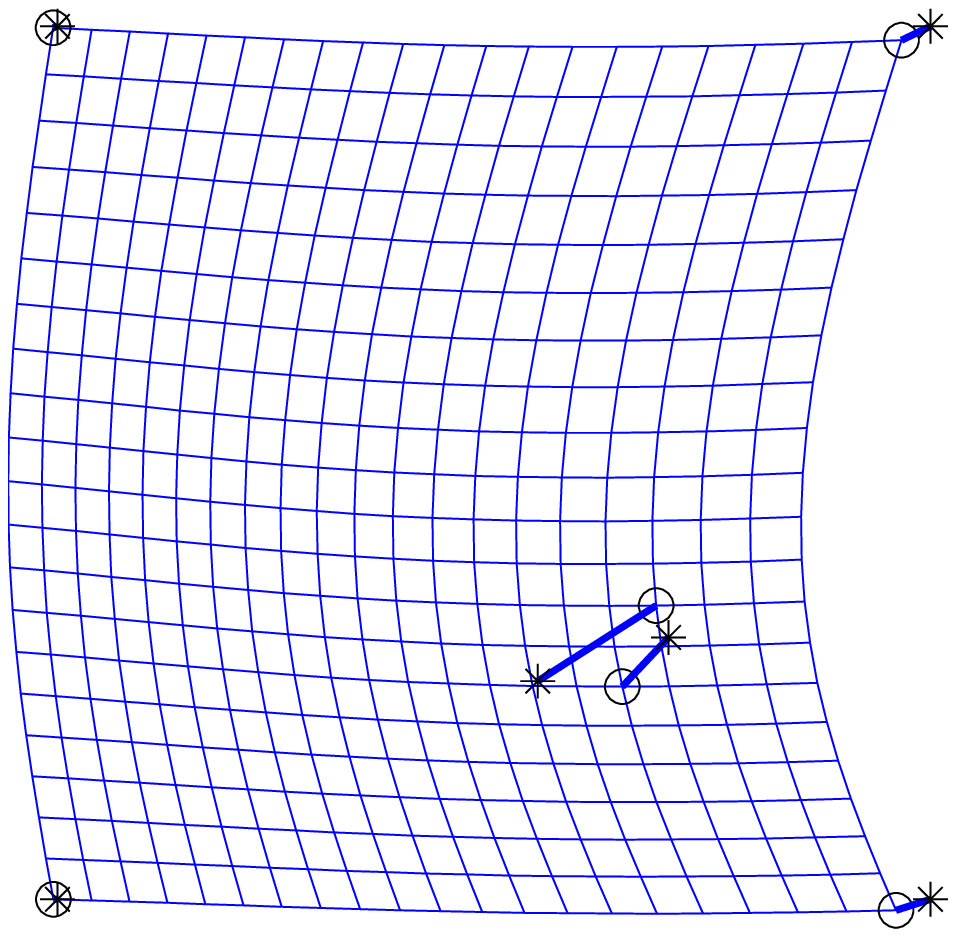} }
\subfigure[] { \includegraphics[width=3cm]{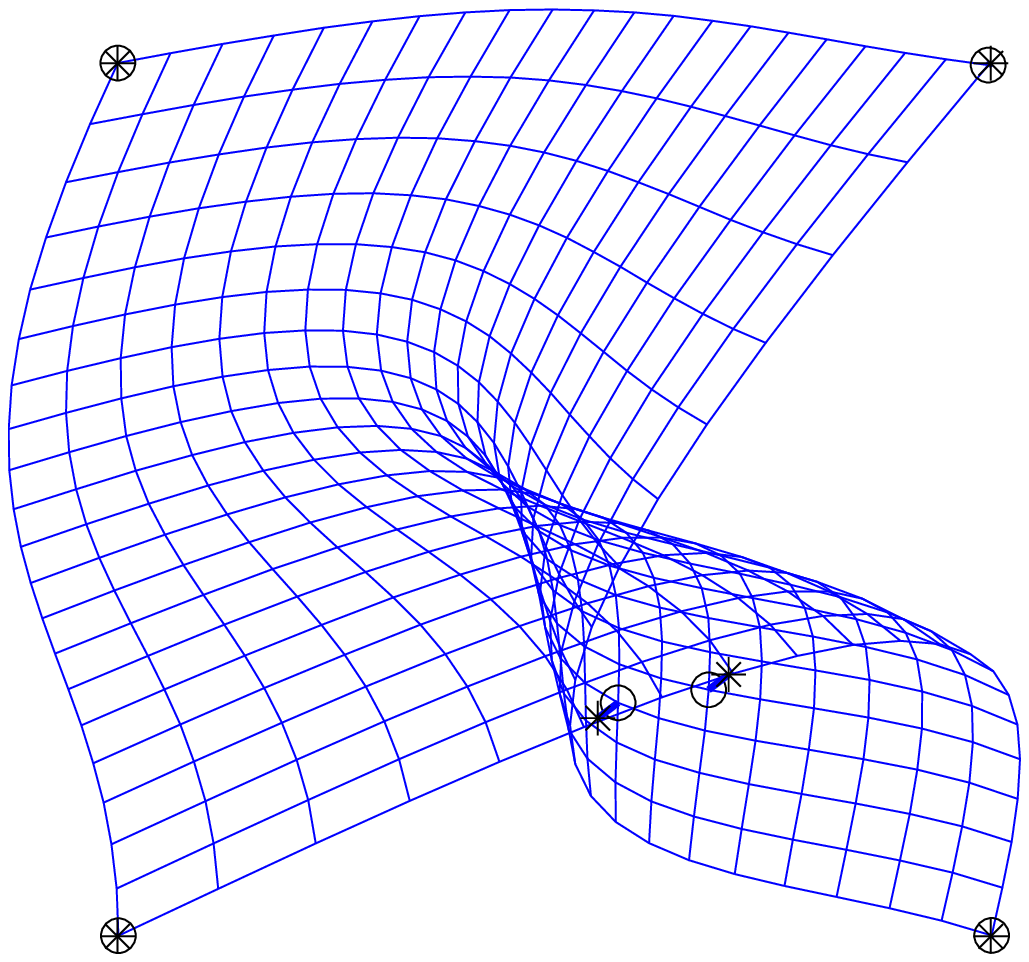} }
\subfigure[] { \includegraphics[width=3cm]{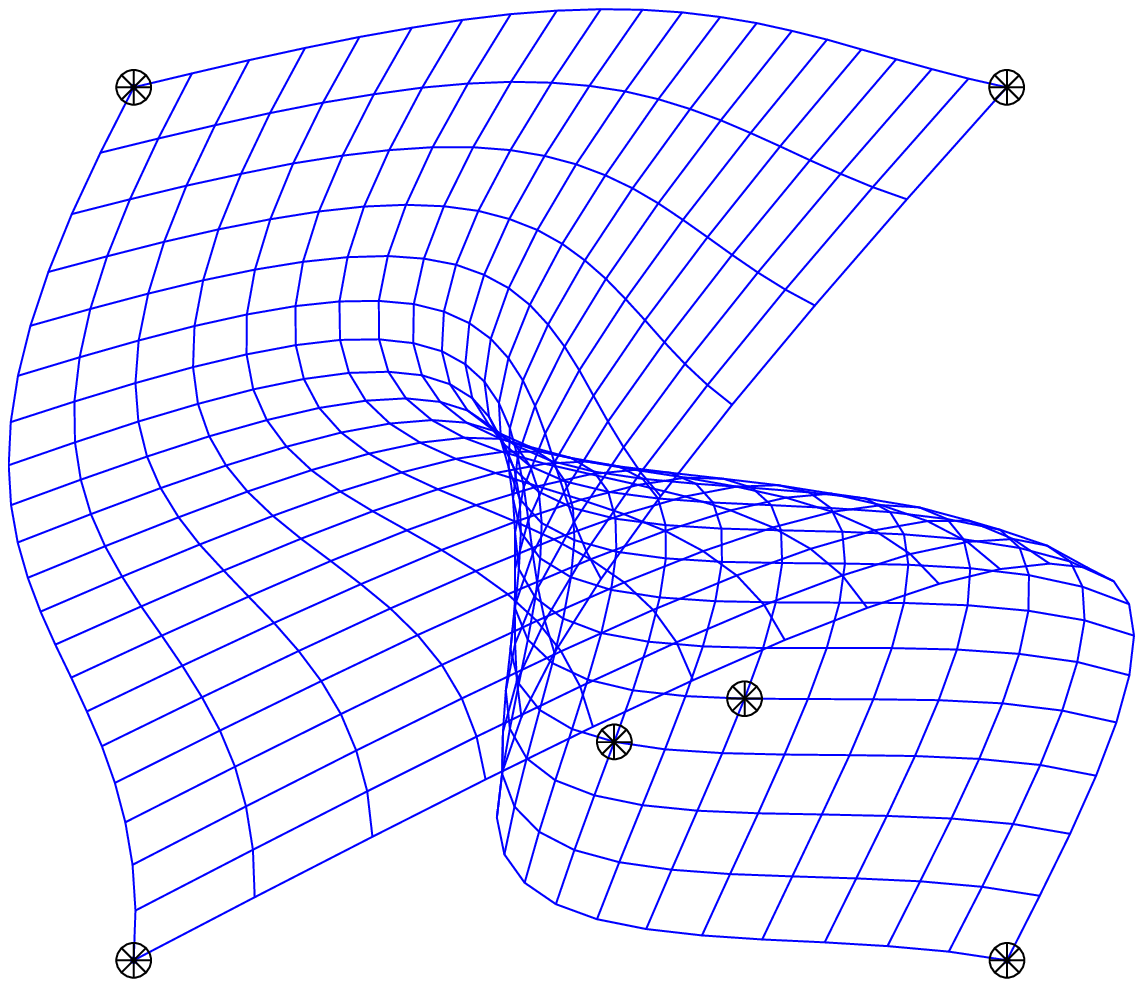} }

\hspace{3cm}
\subfigure[] { \includegraphics[width=3cm]{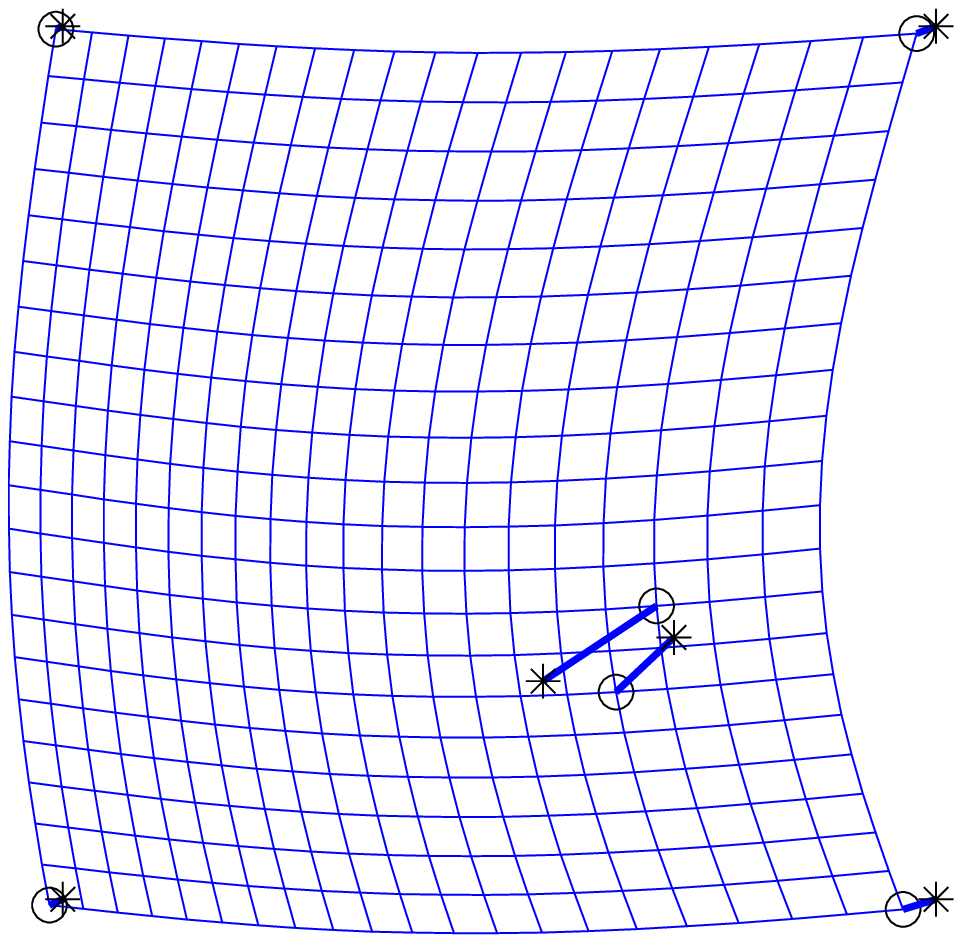} }
\subfigure[] { \includegraphics[width=3cm]{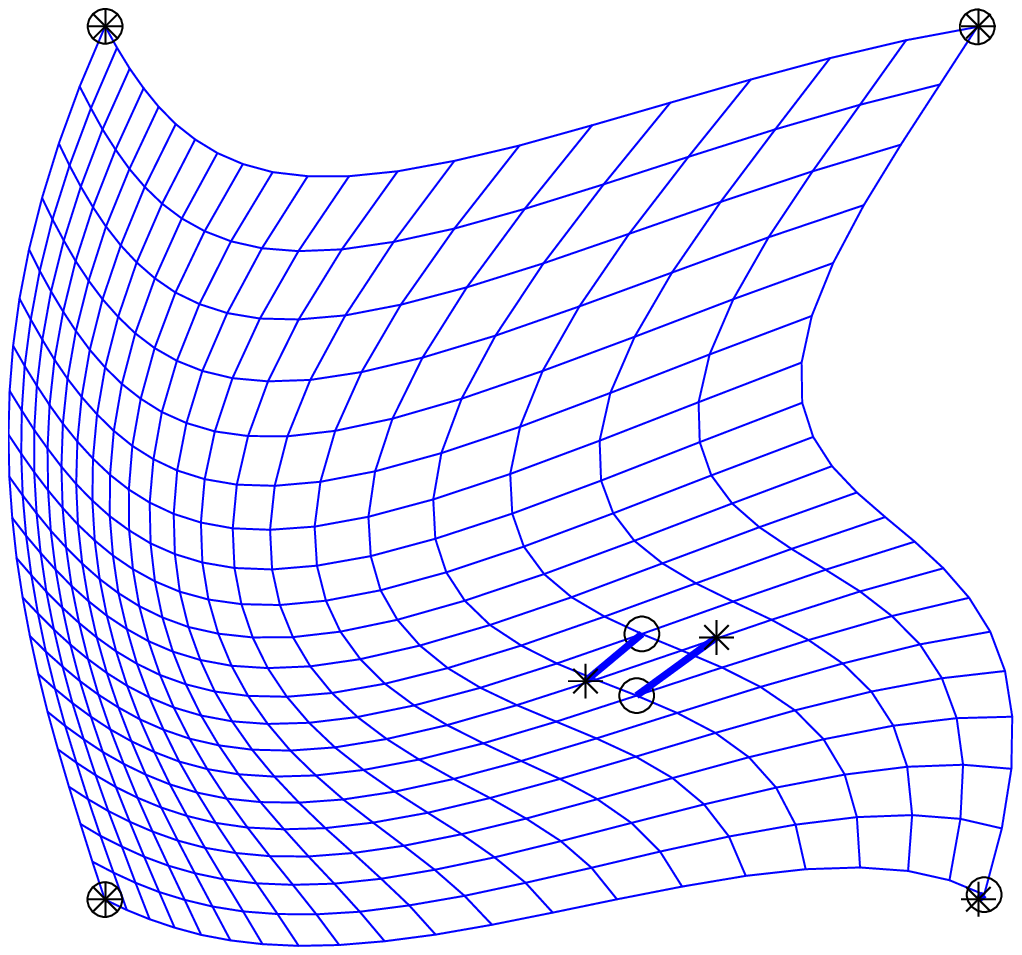} }
\subfigure[] { \includegraphics[width=3cm]{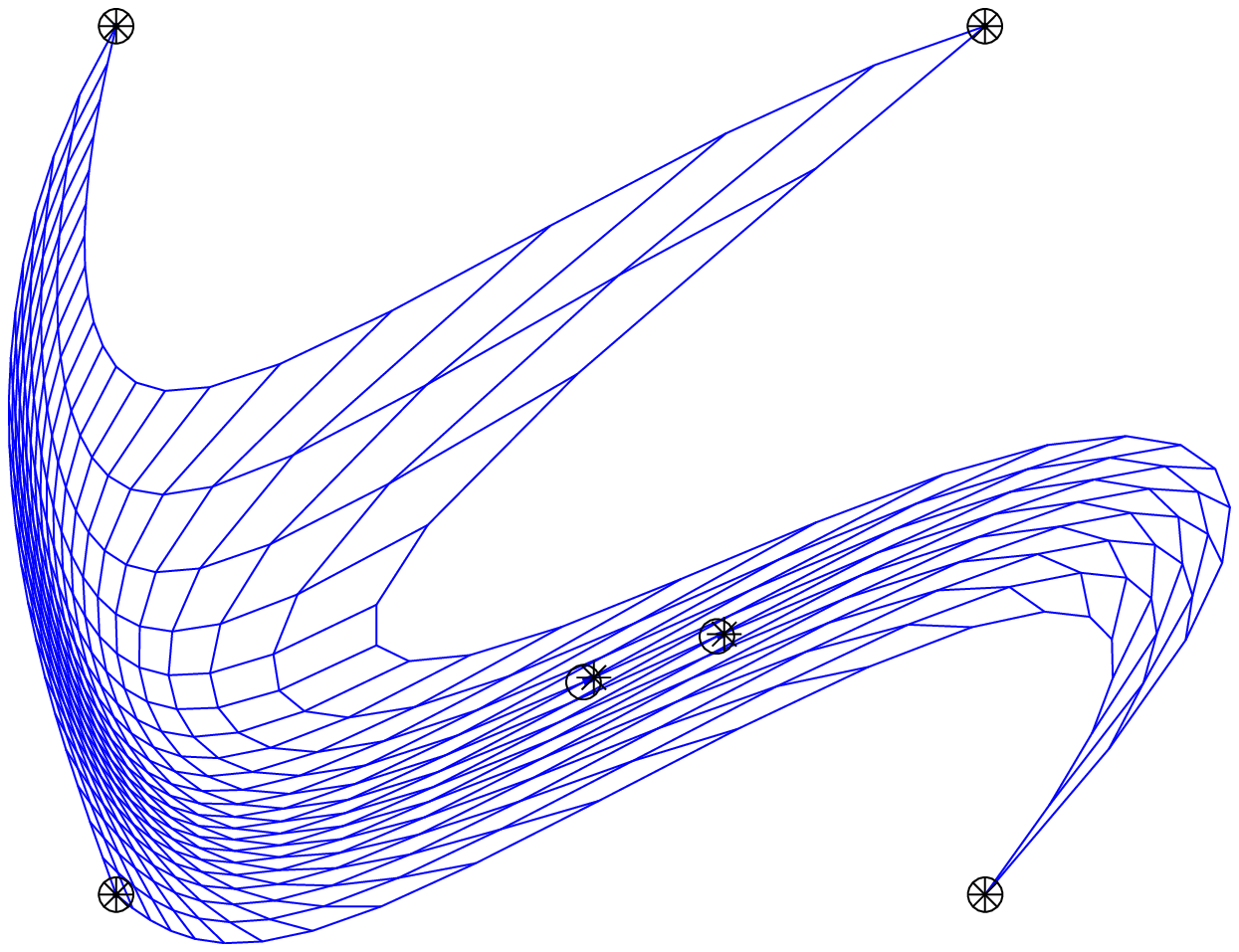} }

\caption{{\em  {\small (a) $n=6$ landmarks to be aligned $x_{1},\ldots,x_{n}$ (circles) and $y_{1},\ldots,y_{n}$ (stars). Four landmarks at the corner of the grid are already at the same location. The bold diagonal lines represents two landmarks be to aligned. Landmark matching with unconstrained spline smoothing:  with a large $\lambda$ (b), a moderate $\lambda$ (c) and a small $\lambda$ (d). Landmark matching with homeomorphic spline:  with a large $\lambda$ (e), a moderate $\lambda$ (f) and a small $\lambda$ (g). }}} \label{fig:alignland}.
\end{figure}

Let $(x_{1},\ldots,x_{n})$ and   $(y_{1},\ldots,y_{n})$ be two sets of $n$ landmark in $\RR^{2}$. The problem of landmark matching is to find a function $f : \RR^{2} \to \RR^{2}$ such that $f(x_{i}) \approx y_{i}$ for all $i=1,\ldots,n$ (see e.g Camion \& Younes \cite{camyou} and references therein). Let $\HH_{K,2}$ be a RKHS of functions in
$L^{2}(\RR^{2})$ with positive definite kernel $K$ and denote by $\tilde{\HH}_{2}$ the set of functions $f : \RR^{2} \to \RR^{2}$ given by
$$
f(x) = Ax + b + \left(\begin{array}{c}  h_{1}(x)  \\  h_{2}(x) \end{array}\right), \quad x \in \RR^{2},
$$
where $h_{1},h_{2} \in  \HH_{K,2} $, $A$ is $2 \times 2$ matrix and $b \in \RR^{2}$.  Landmark matching can be formulated as the problem of finding the minimizer
$$
f_{n,\lambda} = \arg \min_{f \in \tilde{\HH}_{2}} \; \frac{1}{n} \sum_{i=1}^{n} \|f(x_{i})-y_{i}\|^{2}_{\RR^{2}} + \lambda ( \|h_{1} \|^{2}_{K} + \|h_{2} \|^{2}_{K})
$$
Under mild assumptions on the landmarks, the solution of this matching problem is unique and of the form: $\forall x \in \RR, \; f_{n,\lambda}(x) = Ax + b + \left(\begin{array}{c} \sum_{i=1}^{n} \beta_{1,i} K(x,x_{i})   \\  \sum_{i=1}^{n} \beta_{2,i} K(x,x_{i})  \end{array}\right),
$ where $A,b,\beta_1,\beta_2$ are solutions of a simple linear system of
equations (see e.g Camion \& Younes \cite{camyou}). However, there are no constraints in this approach which guarantees that $f_{n,\lambda}$ is a one-to-one mapping of $\RR^{2}$. Indeed, folding are possible for small values of $\lambda$ as shown in Figure \ref{fig:alignland}, where the mapping $f_{n,\lambda}$ is displayed via the deformation of an equally spaced grid of points in $\RR^{2}$.

\subsection{Homeomorphic spline for diffeomorphic matching}

Let $ \tilde{\mathcal{X}}= \{ (v_{t},t \in [0,1]) \mbox{ with }  v_t \in \tilde{\HH}_{2} \mbox{ for all } t \in [0,1]  \}$. Diffeomorphic matching of two sets of landmarks in $\RR^{2}$ by homeomorphic spline is defined as the problem of finding a time-dependent vector field $v \in \tilde{\mathcal{X}}$, which minimizes the ``energy''
\begin{equation}
\tilde{E}_{\lambda}(v) =  \int_{0}^{1} \frac{1}{n} \sum_{i=1}^{n} \left\| y_{i}-x_{i} -v_{t}(ty_{i}+(1-t)x_{i}) \right\|^{2}_{\RR^{2}}dt + \lambda \int_{0}^{1} (\| h_{1,t} \|^{2}_{K}+ \| h_{2,t} \|^{2}_{K}) dt,\label{eqhomeospline2D}
\end{equation}
with $\lambda > 0$ a regularization parameter and
$
v_{t}(x) = A_{t}x + b_{t} + \left(\begin{array}{c}  h_{1,t}(x)  \\  h_{2,t}(x) \end{array}\right),
$
where for each $t \in [0,1]$, $h_{1,t},h_{2,t} \in  \HH_{K,2} $, $A_{t}$ is $2 \times 2$ matrix and $b_{t} \in \RR^{2}$. Then, by taking $ v^{n,\lambda} = \arg \min_{v \in \tilde{\mathcal{X}}} \tilde{E}_{\lambda}(v),$ a diffeomorphic mapping between these two sets of landmarks is obtained by computing
$
f_{n,\lambda}^{c}(x) = \phi^{v^{n,\lambda}}_{1}(x) = x + \int_{0}^{1} v^{n,\lambda}_{t}(\phi^{v^{n,\lambda}}_{t}(x)) dt.
$
Under mild assumptions, and arguing as in the proof of Proposition \ref{propexist}, the optimization problem (\ref{eqhomeospline2D}) has a unique solution $v^{n,\lambda} \in \tilde{\mathcal{X}}$ such that at each time $t \in [0,1]$, $v_{t}^{n,\lambda}$ is the solution of the following standard unconstrained smoothing problem: find $v_{t} \in   \tilde{\HH}_{2}$ which minimizes
$
\tilde{E}^{t}_{\lambda}(v_{t}) = \frac{1}{n}\sum_{i=1}^{n} \left\| y_{i}-x_{i} -v_{t}(ty_{i}+(1-t)x_{i}) \right\|^{2}_{\RR^{2}} + \lambda  (\|h_{1,t} \|^{2}_{K}+\|h_{2,t} \|^{2}_{K}),
$
where $v_{t}(x) =A_{t}x + b_{t} + \left(\begin{array}{c}  h_{1,t}(x)  \\  h_{2,t}(x) \end{array}\right)$. Hence, the computation of the diffeomorphic mapping $f_{n,\lambda}^{c}$ is obtained using unconstrained spline smoothing and by running an ODE. Numerically, we use an Euler scheme similar to the one proposed in Section \ref{seccomp}.

Remark that the formulation (\ref{eqhomeospline2D}) is somewhat similar to the geodesic smoothing spline problem proposed by Camion \& Younes \cite{camyou} in 2D setting. To compute a smooth diffeomorphism to align two sets of landmarks $(x_{i},y_{i}),i=1,\ldots,n$ in $ \RR^{2} \times  \RR^{2}$, Camion \& Younes \cite{camyou} suggest to minimize the following energy
\begin{equation}
\sum_{i=1}^{n} \|\frac{d q_{i}(t)}{dt} -v_{t}(q_{i}(t) \|^{2}_{\RR^{2}} + \lambda  \int_{0}^{1} (\|h_{1,t} \|^{2}_{K}+\|h_{2,t} \|^{2}_{K}) dt \label{eqyou}
\end{equation}
over all time-dependent vector fields and all landmark trajectories $q_{i}(t), i=1,\ldots,n$ with initial conditions $q_{i}(0) = x_{i}$ and $q_{i}(1) = y_{i}$. This leads to an optimization problem which can be solved by a gradient-descent algorithm. In our formulation (\ref{eqhomeospline2D}), the landmarks trajectories are fixed, and correspond to linear paths $q_{i}(t) =  ty_{i}+(1-t)x_{i}$ between $x_{i}$ and $y_{i}$. This makes the optimization problem (\ref{eqyou}) easier to solve.

\begin{figure}[htdp]
\centering \subfigure[] { \includegraphics[width=3cm]{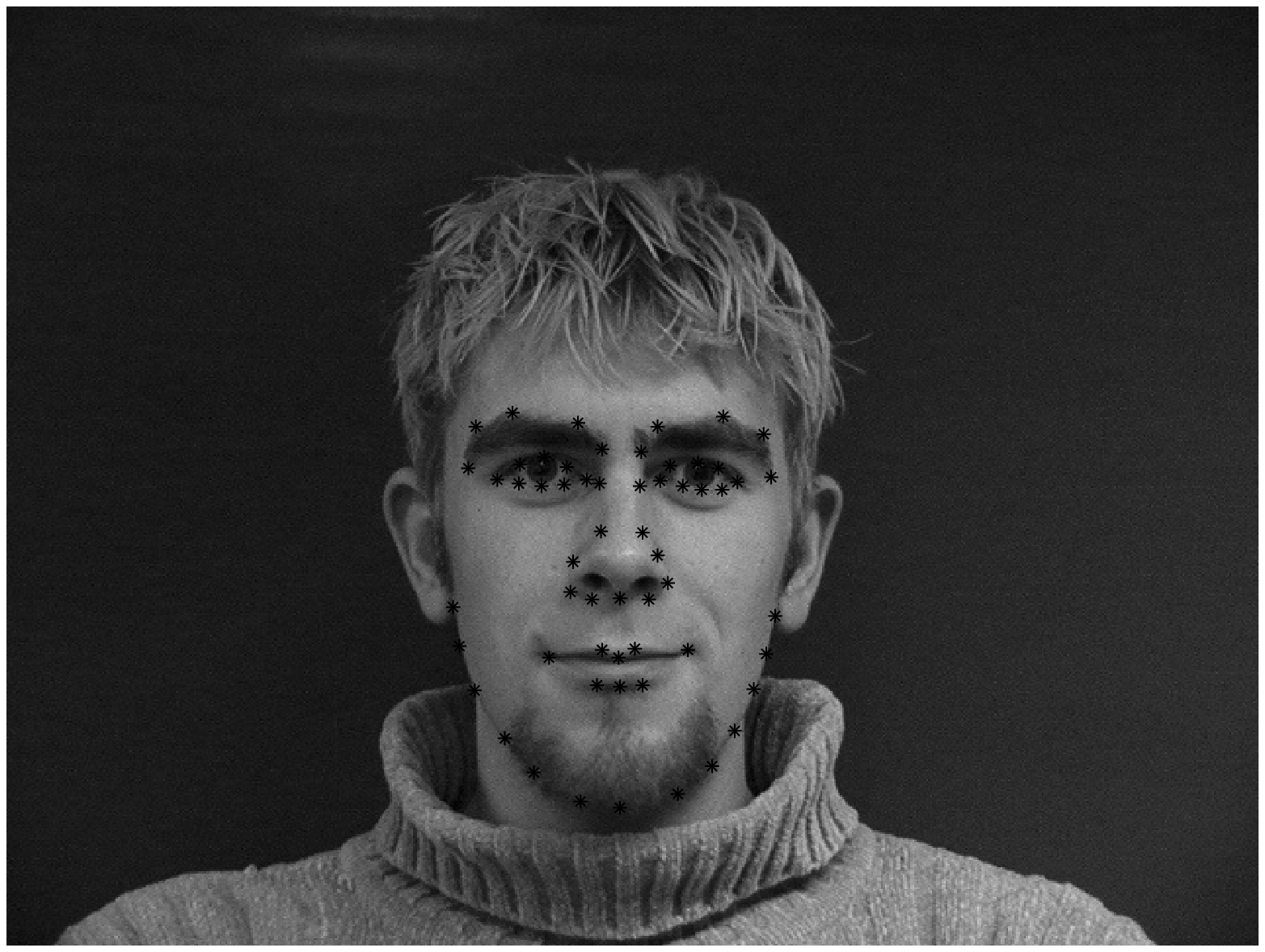} }
\subfigure[] { \includegraphics[width=3cm]{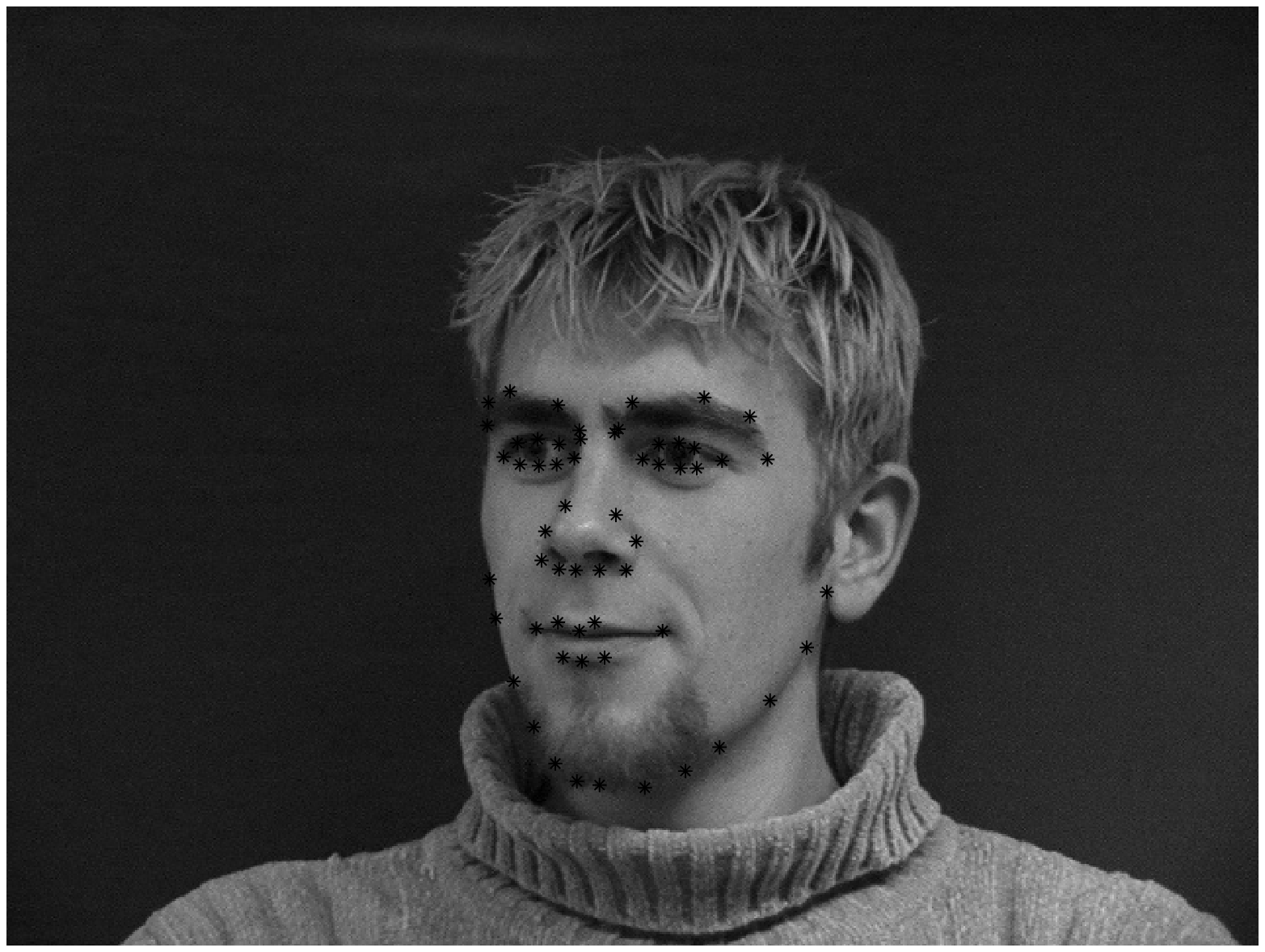} }
\subfigure[] { \includegraphics[width=3cm]{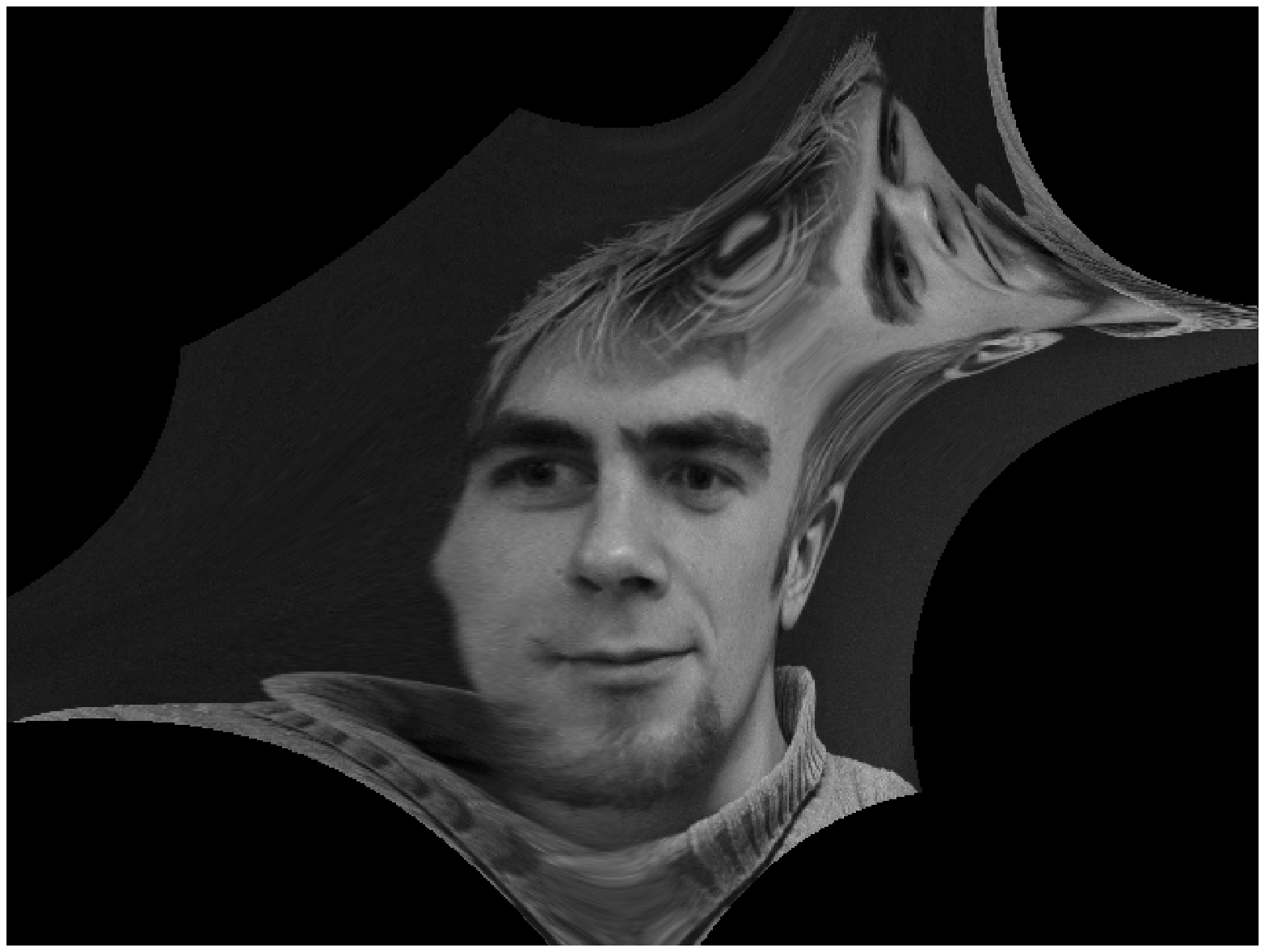} }
\subfigure[] { \includegraphics[width=3cm]{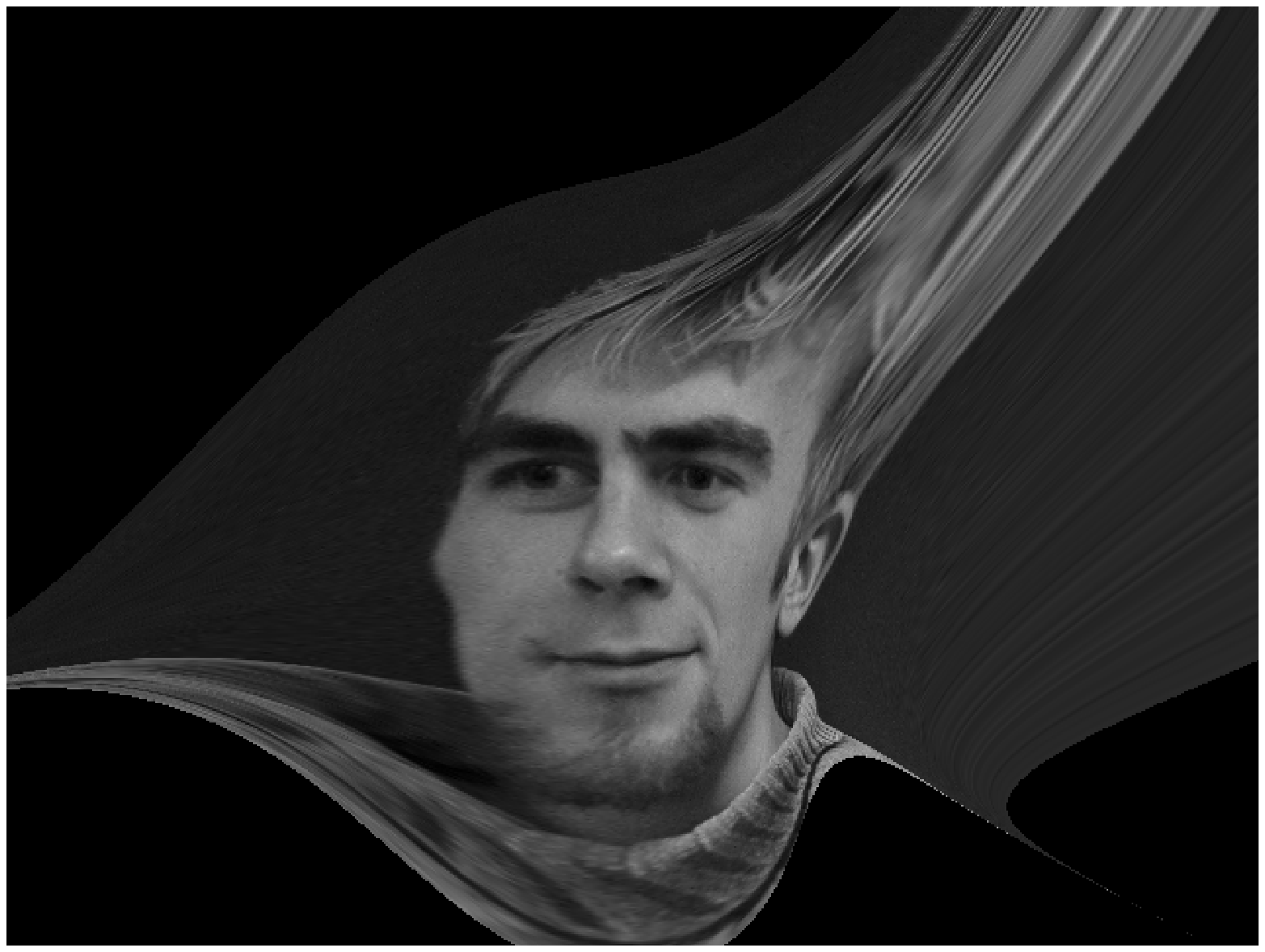} }

\caption{{\em {\small Image warping based on landmark alignment. Figures (a) and (b) are two images with manually annotated landmarks representing facial structures. Figure (c) is the warping of image (b) onto image (a) using landmark alignment with unconstrained spline smoothing. Figure (d) is the warping of image (b) onto image (a) using landmark alignment with homeomorphic spline. Images are taken from the IMM Face Database \cite{steg}.}}} \label{fig:alignimg}
\end{figure}

An example of diffeomorphic mapping is shown in Figure \ref{fig:alignland}, and one can see that even for small values of $\lambda$ the mapping remains one-to-one contrary to unconstrained spline smoothing. An example of landmark-based image warping is also displayed in Figure \ref{fig:alignimg} which illustrates the advantages of homeospline over unconstrained spline smoothing which may lead to unrealistic matching.

\section{Conclusion and Future works}
Homeomorphic splines allows one to compute with a low computational cost, monotone regressors in a 1D setting. In the presence of noisy data, this leads to an estimator which performs at least as well as existing ones. Moreover such an estimator has an optimal rate of convergence over a large class of functional spaces. Homeomorphic splines can also be extended in a 2D setting for landmarks and image warping. In future work, we plan to investigate applications for 2D regression under monotonicity constraints, and also to study the asymptotic normality of our estimator. \\


\noindent {\bf Acknowledgments:} we would like to gratefully acknowledge Christine Thomas-Agnan for discussions about RKHS, and Anestis Antoniadis, Alain Trouv\'e and Laurent Younes for very helpful comments on an early version of the manuscript.  


\section*{Appendix}

\noindent {\bf Proof  of Theorem \ref{thode}:} the key point to derive this result is the use of the following lemma which follows immediately from condition (\ref{eqassLip}) and the assumption $M=2$, $\psi_1(x) = 1$, $\psi_2(x)=x$ made in Section \ref{secRKHS}.


\begin{lemma} \label{lemLipBound}
Suppose that the assumptions of Theorem \ref{thode} are satisfied. Then, for any $g \in \tilde{\HH}$, there exists a
constant $C_{2}$ (not depending on $g$) such that for all $x,y \in \RR$, $|g(x) - g(y)| \leq C_{2} \| g\|  |x-y|$,  and $|g(x)| \leq C_{2} \| g\| (1 + |x|)$.\\
\end{lemma}

\noindent {\bf Existence and uniqueness of the solution of the ODE :} It follows from the Picard-Lindelof Theorem (see Theorem 19 and 20 in Younes \cite{you} using the same notations).

\noindent {\bf Invertibility of the solution :} for $v \in  \mathcal{X}^{1}$, $x \in \Omega$ and any $t,s \in [0,1]$, we denote by $\phi_{st}^{v}(x)$ the value at time $t$ of the unique solution of the equation $\frac{\text{d} \phi_{r}}{\text{d} r} = v_{r}( \phi_{r})$ which is equal to $x$ at time $s$. The solution $\phi_{st}^{v}(x)$ is equal to $\tilde{x} = \phi_{sr}^{v}(x)$ at time $r \in [0,1]$, and thus also equal to $\phi_{rt}^{v}(\tilde{x}) = \phi_{st}^{v}(x)$ which implies that $\phi_{st}^{v}(x) = \phi_{rt}^{v} \circ \phi_{sr}^{v}(x)$ and thus $\phi_{st}^{v} \circ \phi_{ts}^{v}(x) = x$ which proves that for all $t \in [0,1]$ $\phi_{t} = \phi_{0t} $ is invertible, and that its inverse is given by $\phi^{-1}_{t} = \phi_{t0} $.

\noindent {\bf Continuity of $\phi_{t}$ and its inverse: } $\phi_{t}$ is a homeomorphism using Gronwall's lemma and the arguments in Younes  \cite{you}. This ends the proof of Theorem \ref{thode}. $\Box$\\


\noindent {\bf Proof of Theorem \ref{theovf}:} let $t \in [0,1]$ and define the interval
$I_{t} = [tf(0),tf(1)+(1-t)] = [a_{t},b_{t}].$
Given our assumptions on $f$ the function $\phi_{t}(x) = tf(x) + (1-t)x$ is a continuously differentiable homeomorphism from $[0,1]$ to $I_{t}$ such that its inverse is continuous and differentiable with
$
\frac{\text{d} \phi^{-1}_{t}}{\text{d} y}(y) = \frac{1}{tf'(\phi^{-1}_{t}(y))+(1-t)} \mbox{ for all } y \in I_{t},
$
which shows that $\frac{\text{d} \phi^{-1}_{t}}{\text{d} y}$ is continuous since we have assumed that $f'(x) > 0$ for all $x \in [0,1]$. Hence, $\phi_{t}$ is a diffeomorphism from $[0,1]$ to $I_{t}$. Moreover, the above formula for $\frac{\text{d} \phi^{-1}_{t}}{\text{d} y}$ and the fact that $f \in  H^{m}([0,1])$ implies that $\phi^{-1}_{t}$ belongs to $H^{m}(I_{t})$ by application of the chain rule of differentiation.
Let $\tilde{f}_{t}$  the function defined on $I_{t}$ such that
$
\tilde{f}_{t}(y) = f(\phi^{-1}_{t}(y))-\phi^{-1}_{t}(y) \mbox{ for all } y \in I_{t}.
$
Given that  $\phi_{t}$ is a diffeomorphism from $[0,1]$ to $I_{t}$ and that $\phi^{-1}_{t} \in H^{m}(I_{t})$, we can again apply the chain rule of differentiation to show that $\tilde{f}_{t}$  belongs to $H^{m}(I_{t})$. Now, define the sub-space $H_{0}$ of functions in $H^{m}(\RR)$ which coincide with $\tilde{f}_{t}$ on $I_{t} = [a_{t},b_{t}]$. First, we shall prove that $H_{0}$ is not empty. Indeed, choose $c_{t} < a_{t}$ and $d_{t} > b_{t}$ and define the Hermite polynomials $P_{0}$ and $P_{1}$ of degree $2m+1$ such that for all $0 \leq k \leq m-1$, $P_{0}^{(k)}(a_{t}) = \tilde{f}_{t}^{(k)}(a_{t})$ and $P_{0}^{(k)}(c_{t}) = 0$, and $P^{(k)}_{1}(b_{t}) = \tilde{f}_{t}^{(k)}(b_{t})$ and  $P^{(k)}_{1}(d_{t}) = 0$. Then, define the function $v_{t}^{m}$ on $\RR$ such that
$$
v_{t}^{m}(y) = \left\{
\begin{array}{cc}
0 & \mbox{ for } y \in ]-\infty,c_{t}[, \\
P_{0}(y) & \mbox{ for } y \in [c_{t},a_{t}[, \\
\tilde{f}_{t}(y) & \mbox{ for } y \in [a_{t},b_{t}], \\
P_{1}(y) & \mbox{ for } y \in ]b_{t},d_{t}], \\
0 & \mbox{ for } y \in ]d_{t},+\infty[, \\
\end{array}
\right.
$$
By construction of $P_{0}$ and $P_{1}$ and the fact that $\tilde{f}_{t} \in H^{m}(I_{t})$, we have that $v_{t}^{m}$ belongs to $H^{m}(\RR)$, and thus $H_{0}$ is not empty. Since the space $H_{0}$ is closed and convex, it contains a unique element of minimum norm that we denote by  $v^{f}_{t}$ which satisfies  equation (\ref{eqvf1}) by construction of $H_{0}$, which completes the proof Theorem \ref{theovf}.  $\Box$  \vspace{0.2cm}

\noindent {\bf Proof of Proposition \ref{propexist} :} let $t \in [0,1]$. Given our assumptions on the $x_{i}$'s and the $y_{i}$'s we can define the function $v^{n,\lambda}_{t} \in  \tilde{\HH}$ as the smoothing spline which minimizes the energy $E^{t}_{\lambda}(v_{t})$ as defined in equation (\ref{eqener0}). Then, by definition of $v^{n,\lambda}_{t}$ we have that for any $v_{t}  \in   \tilde{\HH}$,  $E^{t}_{\lambda}(v^{n}_{t}) \leq E^{t}_{\lambda}(v_{t})$ which implies that for any $v \in \mathcal{X}$
$
E_{\lambda}(v^{n,\lambda}) = \int_{0}^{1}E^{t}_{\lambda}(v^{n}_{t})dt \leq E_{\lambda}(v),
$
which proves that $v^{n,\lambda}$ is a minimum of  $E_{\lambda}(v)$. Then, the uniqueness of  $v^{n,\lambda}$ follows from the strict convexity of $E_{\lambda}$.

Now, for $i = 1,\ldots,n$ and $t \in [0,1]$, let $\hat{X}_{i}^{t} = ty_{i}+(1-t)x_{i}$ and $\hat{Y}_{i} = y_{i}-x_{i} $. Then, the time-dependent vector field $v_{t}^{n,\lambda}$ is such that for any $x \in \RR$ (see  Wahba \cite{wah})
$v_{t}^{n,\lambda}(x) = \alpha_{1}^{t} + \alpha_{2}^{t}x + \sum_{i=1}^{n} \beta_{i}^{t} K(x,\hat{X}_{i}^{t} ),$
where the coefficients $\balpha_{t} = (\alpha^{t}_{1},\alpha_{2}^{t})'$ and $\bbeta_{t} = (\beta^{t}_{1},\ldots,\beta_{n}^{t})'$ are given by:
$
\bbeta_{t}  = \hat{\Sigma}^{-1}_{\lambda,t}(I_{n}-P_{t}) \bY
$
and
$
\balpha_{t}  = (T_{t}'  \hat{\Sigma}^{-1}_{\lambda,t} T_{t})^{-1}T_{t}'  \hat{\Sigma}^{-1}_{\lambda,t}\bY,
$
with $\bY = (\hat{Y}_{1},\ldots,\hat{Y}_{n})'$, $\hat{\Sigma}_{\lambda,t} = \hat{\Sigma}_{t} + n\lambda I_{n}$ and $\hat{\Sigma}_{t}$ is the $n \times n$ matrix with elements $\hat{\Sigma}_{t}[i,j] = K(\hat{X}_{i}^{t},\hat{X}_{j}^{t})$, $T_{t}$ is the $n \times 2$ matrix with elements $T_{t}[i,1] = 1$ and $T_{t}[i,2] = \hat{X}_{i}^{t}$, and $P_{t}$ is the $n \times n$ matrix given by
$
P_{t} = T_{t}(T_{t}'  \hat{\Sigma}^{-1}_{\lambda,t} T_{t})^{-1}T_{t}'  \hat{\Sigma}^{-1}_{\lambda,t}.
$
Then by the continuity of $K$ and the continuity of matrix inversion, it follows that the coefficients $\balpha_{t}$ and $\bbeta_{t}$ are continuous functions of $t$ on $[0,1]$. Hence, $t \mapsto \|v_{t}^{n,\lambda}\|$ is a continuous function on $[0,1]$ which implies that $\int_{0}^{1} \|v_{t}^{n,\lambda}\|^{2} dt < + \infty$.  Hence  $v^{n,\lambda} \in \mathcal{X}^{2}$ which completes the proof of Proposition \ref{propexist} using Theorem \ref{thode}. $\Box$ \vspace{0.2cm}

\noindent {\bf Proof of Theorem \ref{theomain}:}  in the proof, $C$ will denote a constant whose value may change from line to line.
For $i = 1,\ldots,n$ and $t \in [0,1]$, let $\hat{X}_{i}^{t} = t\hat{f}^{n}(x_{i})+(1-t)x_{i}$,  $X_{i}^{t} = tf(x_{i})+(1-t)x_{i}$ and $\hat{Y}_{i} = \hat{f}^{n}(x_{i})-x_{i} $. One can remark that the smoothing spline $v_{t}^{n,\lambda}$ evaluated at the ``design points'' $\hat{X}_{1}^{t},\ldots,\hat{X}_{n}^{t}$ is a linear function of the observations $\hat{Y}_{1},\ldots,\hat{Y}_{n}$ and can therefore be written as $ \mathbf{v}_{t}^{n,\lambda} = \left(v_{t}^{n,\lambda}(\hat{X}_{1}^{t}),\ldots,v_{t}^{n,\lambda}(\hat{X}_{n}^{t})\right)' = A_{\lambda,t} \bY,$ where $ A_{\lambda,t} = P_{t} + \hat{\Sigma}_{t} \hat{\Sigma}^{-1}_{\lambda,t} (I_{n}-P_{t}).$
Then, under the assumptions of Theorem  \ref{theomain}, the following lemma holds (the proof  follows using arguments in Craven \& Wahba \cite{cra}): 

\begin{lemma} \label{lemeig}
For almost all $t \in [0,1]$ and any $\bx \in \RR^{n}$, we have $ \|A_{\lambda,t} \bx \|_{2} \leq 2 \|\bx \|_{2}$   a.s. 
\end{lemma}


\noindent The next step shows that $v_t^{n,\lambda}$ is in $\mathcal{X}^2$ with an asymptotically probability equals to 1.

\begin{lemma} \label{lemexist}
There exists  $C_{5} > 0$   such that
$
\lim_{n \to + \infty}\PR \left(\int_{0}^{1} \|v_{t}^{n,\lambda}\|^{2} dt \leq   C_{5} \right) = 1.
$
\end{lemma}

\noindent {\bf Proof:} recall that $v^{n,\lambda}$ is the minimizer of the following energy:

$
E_{\lambda_{n}}(v) = \int_{0}^{1}\frac{1}{n}\sum_{i=1}^{n} \left(\hat{f}^{n}(x_{i})-x_{i} -v_{t}(t\hat{f}^{n}(x_{i})+(1-t)x_{i}) \right)^{2}dt + \lambda_{n} \int_{0}^{1} \|h_{t} \|^{2}_{K}dt,
$
where $v_{t}(x) = a_{1}^{t} + a_{2}^{t}x + h_{t}(x)$. We will show that $v^{n,\lambda}$ converges in probability to $v^{f}$ for the norm $\|v\|_{ \mathcal{X}^{2}} = \int_{0}^{1} \|v_{t} \|^{2}dt $. Let $ E(v) =  \int_{0}^{1}  \int_{0}^{1} \left(f(x)-x-v_{t}(t f(x) + (1-t)x) \right)^{2}w(x)dx dt,$
$\mathcal{A}$ be a compact set in $\mathcal{X}^{2}$ and take $v \in \mathcal{A}$:
$
\left( E_{\lambda_{n}}(v) - E(v) \right)^{2} \leq 2E_{1,n}^{2}(v) + 2E_{2,n}^{2}(v),
$
where $E_{2,n}(v) =  \lambda_{n} \int_{0}^{1} \|h_{t} \|^{2}_{K}dt$ and
\begin{eqnarray*}
E_{1,n}(v) &= &  \int_{0}^{1}  \int_{0}^{1} \left(f(x)-x-v_{t}(t f(x) + (1-t)x) \right)^{2}(dw_{n}(x)-w(x)dx) dt +\\
& &  \int_{0}^{1}  \int_{0}^{1} \left[ \left(\hat{f}_{n}(x)-x-v_{t}(t\hat{f}_{n}(x) + (1-t)x \right)^{2} - \right. \\
& & \left. \left(f(x)-x-v_{t}(t f(x) + (1-t)x) \right)^{2} \right]dw_{n}(x) dt,
\end{eqnarray*}
with $w_{n}(x) = \frac{1}{n}\sum_{i=1}^{n} \delta_{x_{i}}(x)$, since $\lambda_{n} \to 0$, we have that
$
\sup_{v \in \mathcal{A}} E_{2,n}^{2}(v) \to 0  \mbox{ as } n \to +\infty. \label{eqE2}
$
Then, remark that $ E_{1,n}^{2}(v)  \leq 2I_{1,n}^{2}(v) + 2I_{2,n}^{2}(v)$, where
$
I_{1,n}(v) = \int_{0}^{1}  \int_{0}^{1} \left(f(x)-x \right.
$
$
\left. -v_{t}(t f(x) + (1-t)x) \right)^{2} (dw_{n}(x)-w(x)dx) dt ,
$
and
$
I_{2,n}(v) =  \int_{0}^{1}  \int_{0}^{1} \left[ \left(\hat{f}_{n}(x)-x \right.\right.
$

$
\left.\left.-v_{t}(t\hat{f}_{n}(x) + (1-t)x ) \right)^{2} -
\left(f(x)-x-v_{t}(t f(x) + (1-t)x) \right)^{2} \right]dw_{n}(x) dt.
$
Let $g^{v}_{t}(x) =  \left(f(x)-x-v_{t}(t f(x) + (1-t)x) \right)^{2}$, then
$
I_{1,n}(v) = \int_{0}^{1}  \int_{0}^{1} g^{v}_{t}(x) dt  (dw_{n}(x)-w(x)dx)
$ by Fubini theorem. Lemma \ref{lemLipBound}, implies that $x  \mapsto \int_{0}^{1} g^{v}_{t}(x) dt$ is bounded on $[0,1]$  by
$
C_{v} = \int_{0}^{1}  \left(1 +\|f\|_{\infty} + C_{2} \|v_{t}\| (2 + \|f\|_{\infty}) \right)^{2} dt.
$
Then, from the compactness of $\mathcal{A}$ it follows that there exists a constant $C$ such that for all $v \in \mathcal{A}$
$
I_{1,n}(v) \leq C \int_{0}^{1} (dw_{n}(x)-w(x)dx).
$
Hence, by definition of $w(x)$ and $w_{n}(x)$ the inequality above finally implies that
\begin{equation}
\sup_{v \in \mathcal{A}} I_{1,n}^{2}(v) \to 0 \mbox{ as } n \to +\infty. \label{eqI12}
\end{equation}
Now, using Cauchy-Schwarz inequality we have that
$
I_{2,n}^{2} \leq I_{3,n}(v)I_{4,n}(v),
$
where
$$
I_{3,n}(v) =  \int_{0}^{1}  \int_{0}^{1} \left(\hat{f}_{n}(x)-f(x)-v_{t}(t\hat{f}_{n}(x) + (1-t)x)+ v_{t}(t f(x) + (1-t)x) \right)^{2}  dw_{n}(x) dt,
$$
and
$
I_{4,n}(v) =  \int_{0}^{1}  \int_{0}^{1} \left(\hat{f}_{n}(x)+f(x)-2x-v_{t}(t\hat{f}_{n}(x) + (1-t)x)- v_{t}(t f(x) + (1-t)x) \right)^{2}  dw_{n}(x) dt.
$
Note that using Lemma \ref{lemLipBound}, it follows that
\begin{equation}
I_{4,n}(v) \leq
2 R_{n}(\hat{f}_{n},f) + 4 \int_{0}^{1} \left(\|f\|_{\infty}+2+2C_{2} \|v_{t}\|(2+\|f\|_{\infty}) \right)^{2}  dt \nonumber    \\
+ 4 C_{2}^{2} R_{n}(\hat{f}_{n},f)  \int_{0}^{1} \|v_{t}\|^{2} dt.
\end{equation}

Then, one has that
$$
I_{3,n}(v) \leq 2  \int_{0}^{1}  \int_{0}^{1} \left(\hat{f}_{n}(x)-f(x)\right)^{2}dw_{n}(x) dt + 2C_{2}^{2}  \int_{0}^{1}  \int_{0}^{1} \|v_{t}\|^{2}\left(\hat{f}_{n}(x)-f(x)\right)^{2} dw_{n}(x) dt,
$$
and thus $ \sup_{v \in \mathcal{A}} I_{3,n}(v)  \leq 2 R_{n}(\hat{f}_{n},f) + 2C_{2}^{2}   \sup_{v \in \mathcal{A}} \int_{0}^{1} \|v_{t}\|^{2}dt R_{n}(\hat{f}_{n},f).$ By assumption, $R_{n}(\hat{f}_{n},f) \to 0$ in probability, and thus $ \sup_{v \in \mathcal{A}} I_{3,n}(v) \to 0 \mbox{ in probability as } n \to + \infty.$ By combining the above equation with the bound for $I_{4,n}(v)$, we finally obtain that  in probability
$
\sup_{v \in \mathcal{A}} I_{2,n}^{2}(v)  \to 0 \mbox{ as } n \to + \infty.
$
So finally, we obtain from (\ref{eqI12})  that $ \sup_{v \in \mathcal{A}} E_{1,n}^{2}(v) \to 0$   in probability  as  $n \to +\infty,$ which implies together with equation (\ref{eqE2}) that $ \sup_{v \in \mathcal{A}} \left( E_{\lambda_{n}}(v) - E(v) \right)^{2}  \to 0 $   in probability as  $n \to +\infty.$
Now, remark that since $E_{\lambda_{n}}(v)$ and $ E(v)$ are positive and strictly convex functionals, they have a unique minimum over the set of time-dependent vector fields $\mathcal{X}$. Moreover, by definition of $v^{f} \in \mathcal{X}^{2}$ one has that for any $t \in [0,1]$, $ v_{t}^{f}(tf(x)+(1-t)x) = f(x)-x$ which implies that $v^{f}$ is the minimum of $E(v)$. Let $\epsilon > 0$ and define $ B(v^{f},\epsilon) = \{v \in \mathcal{X}^{2}\; ; \|v-v^{f}\|_{ \mathcal{X}^{2}} \leq \epsilon\},$ and let $ \text{d} B(v^{f},\epsilon) = \{v \in \mathcal{X}^{2}\; ; \|v-v^{f}\|_{ \mathcal{X}^{2}}  = \epsilon\} $ be the frontier of $B(v^{f},\epsilon)$. Since $v^{f}$ is the minimum of $E(v)$, there exits $\delta > 0$ such that for any $v \in \text{d} B(v^{f},\epsilon)$
$
E(v^{f}) < E(v) + \delta.
$
Obviously, $B(v^{f},\epsilon)$ is a compact subset of $ \mathcal{X}^{2}$, and thus $ \sup_{v \in B(v^{f},\epsilon)} \left|E_{\lambda_{n}}(v) - E(v) \right|$ converges to zero in probability. This implies that for any $\alpha > 0$, there exists $n_{1} \in \mathbb{N}$ such that for any $n \geq n_{1}$,
$
\PR \left( v \in \text{d} B(v^{f},\epsilon)\; ; E_{\lambda_{n}}(v) - E(v) > -\frac{\delta}{3} \right) \geq 1-\frac{\alpha}{2}. 
$
Similarly, there exists $n_{2}$ such that for any $n \geq n_{2}$,
$
\PR \left(  E(v^{f}) - E_{\lambda_{n}}(v^{f}) > -\frac{\delta}{3} \right) \geq 1-\frac{\alpha}{2}.
$
Hence, we have that for any $n \geq \max(n_{1},n_{2})$, $ \PR \left( v \in \text{d} B(v^{f},\epsilon)\; ; E_{\lambda_{n}}(v) > E_{\lambda_{n}}(v^{f}) + \frac{\delta}{3} \right) \geq 1-\alpha.$ This implies that except on a set of probability less than $\alpha$, $E_{\lambda_{n}}$ has a local minimum in the interior of $B(v^{f},\epsilon)$ which is thus the global minimum $v^{n,\lambda}$ since $E_{\lambda_{n}}$ is strictly convex. Hence, we finally have that for any $\epsilon > 0$ then with probability tending to one  $v^{n,\lambda}$ belongs to  $B(v^{f},\epsilon)$ which implies that $v^{n,\lambda}$ converges in probability to $v^{f}$ for the norm $\|\cdot\|_{\mathcal{X}^{2}}$. This proves that there exists a constant $C_{5}$ (not depending on $n$) such that as $n \to + \infty$, $ \PR(\|v^{n,\lambda}\|_{ \mathcal{X}^{2}} \leq C_{5} ) \to 1, $ which completes the proof of Lemma \ref{lemexist}.
$\Box$\\

Now, since $v^{n,\lambda} \in \mathcal{X}^{2}$, Theorem \ref{thode} implies that we can define $\hat{f}_{n,\lambda}^{c}$  and $\phi^{n,\lambda}_{t} $ as the solutions respectively at time $t=1$ and time $t \in [0,1]$ of the ODE $\frac{\text{d} \phi_{t}}{\text{d} t} = v^{n,\lambda}_{t}( \phi_{t})$. 
We shall now  control the empirical error $ R_{n}(\hat{f}_{n}^{c},f) = \frac{1}{n}\sum_{i=1}^{n} (\hat{f}_{n}^{c}(x_{i})-f(x_{i}) )^{2}$. First,
$
\sum_{i=1}^{n} (\hat{f}_{n}^{c}(x_{i})-f(x_{i}) )^{2} =  \sum_{i=1}^{n} (\phi^{n,\lambda}_{1}(x_{i})-\phi_{1}^{f}(x_{i}))^{2}
=  \sum_{i=1}^{n}  \left( \int_{0}^{1}  \left(v^{n,\lambda}_{t}(\phi^{n,\lambda}_{t}(x_{i}))  - v^{f}_{t}(\phi_{t}^{f}(x_{i})) \right)dt \right)^{2}.
$
and note that for any $t \in [0,1]$
$$
\sum_{i=1}^{n} (\phi^{n,\lambda}_{t}(x_{i})-\phi_{t}^{f}(x_{i}))^{2}  =  \sum_{i=1}^{n}  \left( \int_{0}^{t} \left(  v^{n,\lambda}_{s}(\phi^{n,\lambda}_{s}(x_{i}))  - v^{f}_{s}(\phi_{s}^{f}(x_{i})) \right) ds \right)^{2},
$$
which implies that (using Cauchy-Schwarz inequality and the fact that $t \leq 1$)
\begin{eqnarray*}
\sum_{i=1}^{n} (\phi^{n,\lambda}_{t}(x_{i})-\phi_{t}^{f}(x_{i}))^{2}
\leq  2  \sum_{i=1}^{n}   \int_{0}^{t} \left(v^{n,\lambda}_{s}(\hat{X}_{i}^{s}) -  v^{f}_{s}(\phi_{s}^{f}(x_{i})) \right)^{2}ds\\
+ \int_{0}^{t} 4C_{2}^{2}\|v_{s}^{n,\lambda}\|^{2} \sum_{i=1}^{n}  \left(\hat{X}_{i}^{s} - \phi_{s}^{f}(x_{i}) \right)^{2}ds +
\int_{0}^{t} 4C_{2}^{2}\|v_{s}^{n,\lambda}\|^{2} \sum_{i=1}^{n}  \left(\phi^{n,\lambda}_{s}(x_{i})- \phi_{s}^{f}(x_{i}) \right)^{2} ds.
\end{eqnarray*}
To bound this sum, we shall use  the following Lemma whose proof can be found in Younes \cite{you}
\begin{lemma}\label{lemyou}
Consider three continuous and positive functions $c_{s},\gamma_{s}$ and $u_{s}$ defined on $[0,1]$ and such that $ u_{t} \leq c_{t} + \int_{0}^{t} \gamma_{s}u_{s}ds$, then $u_{t} \leq c_{t} + \int_{0}^{t} c_{s} \gamma_{s} e^{\int_{s}^{t}\gamma_{r}dr}ds$.
\end{lemma}

Then, if we apply Lemma \ref{lemyou} by letting
$
u_{t} =  \sum_{i=1}^{n} (\phi^{n,\lambda}_{t}(x_{i})-\phi_{t}^{f}(x_{i}))^{2}, \; \gamma_{t} = 4 C_{2}^{2}\| v^{n,\lambda}_{t}\|^{2},
$
$
c_{t}  =   2  \int_{0}^{t} \sum_{i=1}^{n}   \left(v^{n,\lambda}_{s}(\hat{X}_{i}^{s}) -  v^{f}_{s}(\phi_{s}^{f}(x_{i})) \right)^{2}ds +  \int_{0}^{t} 4C_{2}^{2}\|v_{s}^{n,\lambda}\|^{2} \sum_{i=1}^{n}  \left(\hat{X}_{i}^{s} - \phi_{s}^{f}(x_{i}) \right)^{2}ds,
$
we obtain that
$
u_{t} \leq c_{t} + \int_{0}^{t} c_{s} \gamma_{s} e^{\int_{s}^{t}\gamma_{r}dr}ds.
$
Now recall that by definition of $v^{f}_{s}$ and $\phi_{s}^{f}$, we have $\phi_{s}^{f}(x_{i})  =  X_{i}^{s}$ and $v^{f}_{s}(\phi_{s}^{f}(x_{i}))  =  f(x_{i})-x_{i} =_{def} \tilde{f}(x_{i})$.

Hence,
$
\sum_{i=1}^{n}  \left(v^{n,\lambda}_{s}(\hat{X}_{i}^{s}) -  v^{f}_{s}(\phi_{s}^{f}(x_{i})) \right)^{2} = \|A_{\lambda,s} \bY -  \mathbf{\tilde{f}}  \|_{2}^{2},
$
since by definition of $A_{\lambda,s}$ one has $\left(v_{s}^{n,\lambda}(\hat{X}_{1}^{s}),\ldots,v_{s}^{n,\lambda}(\hat{X}_{n}^{s})\right)' = A_{\lambda,s} \bY$,
and where  $ \mathbf{\tilde{f}} = \left(\tilde{f}(x_{1}),\ldots,\tilde{f}(x_{n}) \right)' $. Similarly, we have that
$
\sum_{i=1}^{n}  \left(\hat{X}_{i}^{s} - \phi_{s}^{f}(x_{i}) \right)^{2} = s^{2}  \|  \hat{\bF}_{n} -  \bF \|_{2}^{2} \leq \|  \hat{\bF}_{n} -  \bF \|_{2}^{2},  
$
where $\hat{\bF}_{n} = \left(\hat{f}_{n}(x_{1}),\ldots,\hat{f}_{n}(x_{n}) \right)' $ and $\bF = \left(f(x_{1}),\ldots,f(x_{n}) \right)' $. Hence,
\begin{equation}
c_{t} \leq 2  \int_{0}^{t}  \|A_{\lambda,s} \bY -  \mathbf{\tilde{f}}  \|_{2}^{2} ds + \int_{0}^{t} 4C_{2}^{2}\|v_{s}^{n,\lambda}\|^{2}ds  \|  \hat{\bF}_{n} -  \bF \|_{2}^{2} . \label{eqct3}
\end{equation}
Now, remark that
$
\|A_{\lambda,s} \bY -  \mathbf{\tilde{f}}  \|_{2}^{2}  \leq \|(I_{n}-A_{\lambda,s}) \mathbf{\tilde{f}}  \|_{2}^{2} + \|A_{\lambda,s}( \bY -  \mathbf{\tilde{f}}) \|_{2}^{2},
$
and observe that by Lemma \ref{lemeig}
\begin{equation}
\|A_{\lambda,s}( \bY -  \mathbf{\tilde{f}})  \|_{2}^{2} \leq 2 \|  \hat{\bF}_{n} -  \bF \|_{2}^{2}, \label{eqI1a}
\end{equation}
and let $\tilde{v}_{s}^{n,\lambda}$ be the solution of the following smoothing problem: find $v \in   \tilde{\HH}$ which minimizes
$
\frac{1}{n}\sum_{i=1}^{n} \left(\tilde{f}(x_{i}) -v(\hat{X}_{i}^{s} ) \right)^{2} + \lambda  \|h\|^{2}_{K}, 
$
where $v_{t}(x) = a_{1}^{t} + a_{2}^{t}x + h_{t}(x)$. Then, by definition of $\tilde{v}_{s}^{n,\lambda}$ we have that 
$$
\frac{1}{n} \|(I_{n}-A_{\lambda,s}) \mathbf{\tilde{f}}  \|_{2}^{2} + \lambda  \| \tilde{h}_{s}^{n,\lambda} \|^{2}_{K}  \leq  \frac{1}{n} 2 C_{2}^{2} \|v_{s}^{f}\|^{2} \sum_{i=1}^{n}  \left(X_{i}^{t} -\hat{X}_{i}^{t} \right)^{2} + \lambda  \| h_{s}^{f} \|^{2}_{K}
$$
Finally, by combing equations (\ref{eqI1a}) and  the above inequality, we obtain that there exists a constant $C$ such that
\begin{equation}
\|A_{\lambda,s} \bY -  \mathbf{\tilde{f}}  \|_{2}^{2}    \leq  C \left( \|  \hat{\bF}_{n} -  \bF \|_{2}^{2} +   n \lambda \right) \mbox{ with } C = \max( 2 C_{2}^{2} \|v_{s}^{f}\|^{2},\| h_{s}^{f} \|^{2}_{K}). \label{eqA}
\end{equation}
Now using Lemma \ref{lemexist}, and combining the above relation with equations (\ref{eqct3}) and (\ref{eqA}), we finally obtain in probability
$
c_{t} \leq b(n).
$
where $b(n) =  C\left( \|  \hat{\bF}_{n} -  \bF \|_{2}^{2} +   n \lambda \right)$ for some constant $C > 0$. Similarly, by Lemma \ref{lemexist} we have that there exists a constant $C_{6}$ such that (in probability) for any $s,t \in [0,1]$
$
\int_{s}^{t}\gamma_{r}dr \leq C_{6}.
$
Then, by combining the previous inequalities, we derive from Lemma (\ref{lemyou}) that $ u_{t} \leq b(n) (1 + C_{6}e^{ C_{6}}).$ Now, since $R_{n}(\hat{f}_{n}^{c},f) = \frac{1}{n} u_{1}$ and $R_{n}(\hat{f}_{n},f) = \frac{1}{n} \|  \hat{\bF}_{n} -  \bF \|_{2}^{2}$, we finally obtain that there exists a constant $\Lambda_{1}$ such that with probability tending to one as $n \to +\infty$
$
R_{n}(\hat{f}_{n}^{c},f) \leq \Lambda_{1} \left( R_{n}(\hat{f}_{n},f) +    \lambda_{n} \right),
$
which completes the proof of Theorem \ref{theomain}. $\Box$

\end{document}